\documentclass[a4paper,10pt]{amsart}
\usepackage[english]{babel}
\usepackage{graphicx}
\usepackage{amsmath}
\usepackage{amssymb}
\usepackage{amsthm}
\usepackage{libertine}
\usepackage[all]{xy}
\usepackage{bm}
\usepackage{anysize}
\usepackage{graphicx}
\usepackage{color}
\usepackage{fancybox}
\usepackage{bbm}
\usepackage{amsfonts}
\usepackage[T1]{fontenc}
\usepackage{latexsym}
\usepackage{cite}
\usepackage[numbers]{natbib}
\usepackage{mathrsfs} 

\usepackage{pst-node}
\usepackage{tikz-cd}
\usepackage{tikz}
\usepackage{float}

 \newtheorem{thm}{Theorem}[section]
 \newtheorem{cor}[thm]{Corollary}
 \newtheorem{lem}[thm]{Lemma}
 \newtheorem{prop}[thm]{Proposition}
  
 \theoremstyle{definition}
 \newtheorem{defi}[thm]{Definition}
 \theoremstyle{remark}
 \newtheorem{rem}[thm]{Remark}
 
 \numberwithin{equation}{section}

\newcommand{\N}{\mathbb{N}}

\newcommand{\R}{\mathbb{R}}      
\newcommand{\Z}{\mathbb{Z}}
\newcommand{\C}{\mathbb{C}}

\newcommand*{\m}[1]{\underline{#1}}

\newcommand{\fd}{\rightarrow}

\newcommand{\inc}{\subset}

\newcommand{\iso}{\cong}

\newcommand{\ba}{\overline}

\DeclareMathAlphabet{\mathpzc}{OT1}{pzc}{m}{it}
\newcommand{\al}{\alpha}
\newcommand{\be}{\beta}
\newcommand{\lan}{\lambda}

\newcommand{\Del}{\Delta}
\newcommand{\gam}{\gamma}
\newcommand{\Gam}{\Gamma}

\newcommand{\Om}{\Omega}

\newcommand{\Ha}{\mathbb{H}}

\newcommand{\Sa}{\mathbb{S}}

\newcommand{\pa}{\partial}

\newcommand{\sgn}{\mbox{sgn}}

\newcommand{\el}{\ell}

\title{On the connection between the Fueter-Sce-Qian theorem and the generalized CK-extension}
\author[De Martino]{Antonino De Martino}
\address{Dipartimento di Matematica,  Politecnico De Milano, Via Bonardi 9, 20133, Milan, Italy}
\email{antonino.demartino@polimi.it}
\author[Diki]{Kamal Diki}
\address{Schmid College of Science and Technology, Chapman University, One University Drive, Orange 92867, California, USA}
\email{diki@chapman.edu}
\author[Guzmán Adán]{Alí Guzmán Adán*}
\thanks{*Postdoctoral Fellow of the Research Foundation - Flanders (FWO)}
\address{Clifford Research Group, Department of Electronics and Information Systems,
 Faculty of Engineering and Architecture, Ghent University, Krijgslaan 281, 9000 Gent, Belgium.}
\email{ali.guzmanadan@ugent.be}

\begin{document}
\begin{abstract}
The Fueter-Sce-Qian theorem provides a way of inducing axial monogenic functions in $\mathbb{R}^{m+1}$ from holomorphic intrinsic functions of one complex variable. This result was initially proved by Fueter and Sce for the cases where the dimension $m$ is odd using pointwise differentiation, while the extension to the cases where $m$ is even was proved by Qian using the corresponding Fourier multipliers. 

In this paper, we provide an alternative description of the Fueter-Sce-Qian theorem in terms of the generalized CK-extension. The latter characterizes axial null solutions of the Cauchy-Riemann operator in $\mathbb{R}^{m+1}$ in terms of their restrictions to the real line. This leads to a one-to-one correspondence between the space of axially monogenic functions in $\mathbb{R}^{m+1}$ and the space of analytic functions of one real variable. 

We provide explicit expressions for the Fueter-Sce-Qian map in terms of the generalized CK-extension for both cases, $m$ even and $m$ odd. These expressions allow for a plane wave decomposition of the Fueter-Sce-Qian map or, more in particular, a factorization of this mapping in terms of the dual Radon transform. In turn, this decomposition provides a new possibility for extending  the Coherent State Transform (CST) to Clifford Analysis. In particular, we construct an axial CST defined through the Fueter-Sce-Qian mapping, and show how it is related to the axial and slice CST’s already studied in the literature.

\vspace{0.3cm}


\end{abstract}

\maketitle

\tableofcontents

\section{Introduction}

The aim of this paper is to bring together two essential results in Clifford analysis: the Fueter theorem   \cite{MR1509515, MR97386, MR1485323, MR4240465} and the generalized Cauchy–Kovalevskaya extension  \cite{MR3077647, MR1169463, CK_Ali} (CK-extension for short). There is a vast literature on both results, which essentially provide two different ways of transforming analytic functions of one (real or complex) variable into monogenic functions, i.e.\ null-solutions of the generalized Cauchy-Riemann operator $\mathcal{D}_x$ in $\R^{m+1}$. The later is a conformally invariant first-order differential operator which both extend the complex Cauchy-Riemann operator $\pa_{\ba{z}}$ to higher dimensions and also factorizes the Laplace operator in $\R^{m+1}$. In explicit coordinates, this operator is defined as $\mathcal{D}_x := \pa_{x_0} + \sum_{j=1}^m e_j \pa_{x_j}$, where $\{e_1, \ldots, e_m\}$ is an orthonormal basis in $\R^m$ underlying the construction of the real Clifford algebra $\R_m$. See e.g.\ \cite{MR697564, MR1169463} for standard references in the Clifford analysis setting.

In its original form  \cite{MR1509515}, Fueter's theorem asserts that every intrinsic holomorphic function induces a quaternionic monogenic function. We recall that a holomorphic function $f(z) = \al(u,v) + i\be(u, v)$ ($z=u+iv$) is said to be intrinsic if it is defined on a complex domain symmetric with respect to the real axis, and if its restriction $f\large|_\R$ to the real line is $\R$-valued (a more detailed definition is provided in Definition \ref{IntHolDef}). Under these conditions, Fueter's result states that the function
\begin{equation}\label{QuatFuet}
F(q_0, \m{q}) = \Del_q \left(\al(q_0,|\m{q}|) + \frac{\m{q}}{|q|}\be(q_0, |\m{q}|)\right)
\end{equation}
is monogenic with respect to the quaternionic Cauchy-Riemann operator 
\[D_q = \pa_{q_0} + {\bf i} \pa_{q_1}+{\bf j}\pa_{q_2}+ {\bf k}\pa_{q_3},\]
 i.e.\ $D_q F(q_0, q) =0$ . Here $\Del_q= \pa_{q_0}^2 + \pa_{q_1}^2+\pa_{q_2}^2+ \pa_{q_3}^2$ is the four-dimensional Laplacian, and $q =  {q_0} + {\bf i} {q_1}+{\bf j}{q_2}+ {\bf k}{q_3}$ is a quaternionic variable whose real and imaginary parts are given by $q_0$ and $\m{q}= {\bf i} {q_1}+{\bf j}{q_2}+ {\bf k}{q_3}$ respectively.
 
 This result was later generalized to $\R^{m+1}$ by Sce \cite{MR97386, MR4240465} for odd values of the dimension $m$. Sce's extension asserts that, given a holomorphic complex function $f$ as above, the $\R_m$-valued function 
\begin{equation}\label{Fuet-Sce-R}
G(x_0, \m{x}) = \Del^{\frac{m-1}{2}} \left(\al(x_0,|\m{x}|) + \frac{\m{x}}{|x|}\be(x_0, |\m{x}|)\right), 
\end{equation}
is in the kernel of the Cauchy-Riemann operator $\mathcal{D}_x$ in $\R^{m+1}$ . Here $\Del= \sum_{j=0}^m \pa_{x_j}^2$ is the Laplacian in $m+1$ real variables, and $x=(x_0,\m{x})$ is a vector variable in $\R^{m+1}$ with $\m{x}= (x_1, \ldots, x_m)\in \R^m$. 


In \cite{MR1485323}, Qian extended Fueter's and Sce's results to all dimensions $m\in\N$ by using the Fourier multiplier definition of the fractional Laplacian for the cases where $m$ is even (which is also valid for the cases where $m$ is odd). This approach leads to the so-called  Fueter-Sce-Qian Theorem (see Theorems \ref{FueterThm1} - \ref{FSQThm}) which establishes a link between holomorphic intrinsic functions and monogenic functions in $\R^{m+1}$ regardless of the parity of $m$. Instead of using the pointwise differential operator  $\Del^{\frac{m-1}{2}}$, this link is established by means of the Fueter-Sce-Qian mapping $\tau_m$  which coincides with $\Del^{\frac{m-1}{2}}$ when $m$ is odd (see Section \ref{SecFuet} 
for a precise definition).

There are two important facts about this result that will be crucial in the following. First of all, the range of the Fueter-Sce-Qian mapping is proven to be the module of monogenic functions of axial type, i.e. functions of the form $A(x_0,|\m{x}|) + \frac{\m{x}}{|\m{x}|}B(x_0, |\m{x}|)$ where $A,B$ are $\R_m$-valued functions. We shall denote the  $\R_m$-module of axial monogenic functions in $\R^{m+1}$ by $\mathcal{AM}(\R^{m+1})$ (see Definition \ref{AxMonDef}).

The second important fact about the Fueter-Sce-Qian extension is that it is realized in two steps. First, one converts the holomorphic complex function $f$ to a Clifford-valued function defined in $\R^{m+1}$ by means of the map 
\begin{equation}\label{SExt1}
 \al(u,v) + i\be(u, v) \;\;\mapsto \;\; \al(x_0,|\m{x}|) + \frac{\m{x}}{|x|}\be(x_0, |\m{x}|), 
\end{equation}
 and later, one acts with the mapping $\tau_m$.
 

The map in (\ref{SExt1}) is known as the slice extension map since it maps intrinsic holomorphic functions in $\C$ to the $\R_m$-module $\mathcal{SM}(\R^{m+1})$ of so-called slice monogenic functions in $\R^{m+1}$ (see Section \ref{SMFSect}). Observe that this map can naturally be defined on $\R$-valued functions of one real variable. This is  due to the equivalence between the space of intrinsic holomorphic functions in $\C$, and the space $\mathcal{A}(\R)$ of $\R$-valued real-analytic functions on $\R$ that can be holomorphically extended to 
the entire complex plane $\C$. In this way, we obtain a map $S: \mathcal{A}(\R)\otimes \R_m \fd  \mathcal{SM}(\R^{m+1})$  which factorizes the Fueter-Sce-Qian extension as follows
%
\begin{equation*}
\begin{tikzcd}[row sep = 3em, column sep = 6em]
\mathcal{A}(\R)\otimes \R_m \arrow[r, "S", rightarrow]& \mathcal{SM}(\R^{m+1}) \arrow[d, "\tau_m"] \\
              & \mathcal{AM}(\R^{m+1}) 
\end{tikzcd}
\end{equation*}

For other generalizations of Fueter's theorem beyond Sce and Qian we refer the reader to \cite{MR3265285, DixanThesis, MR1957490, MR2234105, MR1812650}. 

The other monogenic extension map studied in this paper is the so-called generalized CK-extension  \cite{MR1169463, MR3077647}. This result 
asserts that (in a similar way to intrinsic holomorphic functions) axially monogenic functions are completely determined by their restrictions to the real line. Conversely, any real analytic function $f(x_0) \in \mathcal{A}(\R)$ has a unique axial monogenic extension GCK$[f](x_0,\m{x})$ to $\R^{m+1}$ (see Theorem \ref{GCKthm}). This leads to an isomorphism  between modules $\textup{GCK} : \mathcal{A}(\R)\otimes \R_m \fd \mathcal{AM}(\R^{m+1})$, which can be included in the previous diagram as follows 
\begin{equation}
\begin{tikzcd}[row sep = 3em, column sep = 6em]
\mathcal{A}(\R)\otimes \R_m \arrow[r, "S", rightarrow] & \mathcal{SM}(\R^{m+1}) \arrow[d, "\tau_m"] \\
\mathcal{A}(\R)\otimes \R_m \arrow[r, "\textup{GCK}", rightarrow]                 & \mathcal{AM}(\R^{m+1})
\end{tikzcd} .
\end{equation}
It is thus natural to try to relate  the mappings GCK and $\tau_m\circ S$. The first goal of this paper is to completely answer the question of  how we can complete and make the above diagram commutative. As we will show throughout Section \ref{FuetGCKSect}, the solution to this problem consists in drawing a vertical downwards arrow in the left side corresponding (up to a constant) to the $(m-1)$-th derivative 
$\pa_{x_0}^{m-1}: \mathcal{A}(\R)\otimes \R_m \fd \mathcal{A}(\R)\otimes \R_m$. 
Strictly speaking, we will show that 
\[\tau_m\circ S = \gam_m \, \textup{GCK} \circ \pa_{x_0}^{m-1},\] 
where $\gam_m \in \C$ is a suitable constant coefficient. 

It is worth mentioning that this idea has been implicitly used for computing the Fueter primitives of the kernels involved in the Fueter's inversion integral formula, see e.g.\ \cite{MR3539492, MR2787441, MR3047080}. In this manuscript, however, we provide a more general and systematic approach for this procedure (see Theorems \ref{res1} and \ref{GenRel}). Moreover, we use our approach to further study the Fueter-Sce-Qian mapping (see Sections \ref{RadTSect} and \ref{CSTSect}). Indeed, as a direct application of our results, we are able to provide explicit expressions for the actions of the Fueter-Sce-Qian mapping on monomials of the form $(x_0+\m{x})^k$, $k\in\N$. This is done by showing that the polynomials $\tau_m\left[(x_0+\m{x})^k\right]$ coincide (up to constant coefficients) with a well-known family of monogenic Appell polynomials introduced in \cite{MR2769480} (see Section \ref{APSect}). {In the forthcoming paper \cite{F-S_Ant_Kam_Ali}, we further apply these results to study the extensions of the classical Fock and Hardy spaces to the Clifford analysis setting.}

Another important application of our results is related to the dual Radon transform $\check{R}$ \cite{MR754767}.  In \cite{MR3412341}, it was shown that the dual Radon transform maps slice monogenic functions into axial monogenic functions, i.e.
\[
\check{R}:\mathcal{SM}(\R^{m+1}) \fd \mathcal{AM}(\R^{m+1}).
\]
Thus the natural question arises as whether it is possible to relate the 
dual Radon transform $\check{R}$ to the Fueter-Sce-Qian mapping $\tau_m$. Answering this question is the second goal of the present paper. Throughout Section \ref{RadTSect}, we shall show that such a relation can be derived from a plane wave decomposition of the generalized CK-extension (see Theorem \ref{PWDGCK}). This leads to a plane wave decomposition of the Fueter-Sce-Qian mapping, which in turn provides a full overview on how the operators GCK, $S$, $\tau_m$ and $\check{R}$ are related to each other. Such relations are summarized in Theorem \ref{FullRelRadThm} and depicted in the diagram (\ref{DiagComp}).

The final goal of this paper is motivated by a question left open in \cite{MR3556034} about the possibility of constructing Segal-Bargmann or Coherent States Transforms (CSTs for short) by means of the  Fueter-Sce-Qian mapping. We recall that the Segal-Bargmann transform on $\R$  \cite{MR538026, MR143519, MR201959} establishes a unitary isomorphism between the Hilber space $\mathcal{L}^2(\R,dx_0)$ and the space $\mathcal{HL}^2 (\C, e^{-y^2}dx_0dy)$, of entire holomorphic functions in $\C$ that are square integrable with respect to the measure $ e^{-y^2}dx_0dy$.
 This transform corresponds to applying first a convolution with the heat kernel in $\R$, followed by a holomorphic extension from $\R$ to $\C$. By adapting this methodology, new CSTs have been defined of a wide variety of manifolds such as compact connected Lie groups, see e.g.\  \cite{MR1274586}. 

In \cite{MR3556034}, two extensions, $U_s$ and $U_a$, of the Segal-Bargmann transform were introduced from $\mathcal{L}^2(\R,dx_0)\otimes \R_m$ to the modules of slice monogenic and axial monogenic functions respectively. These transforms were defined by means of the slice and the GCK extension maps, and are related by means the dual Radon transform, i.e.
\[
U_a = \check{R} \circ U_s.
\] 
In \cite[Remark 4.6]{MR3556034}, the authors proposed an interesting alternative for defining yet another axial CST through the Fueter-Sce-Qian theorem, i.e. a CST using the  Fueter-Sce-Qian mapping as extension map from $\R$ to $\R^{m+1}$ instead of the GCK extension. The study of a possible relation of such a transform with the ones already defined in the monogenic setting  was left as an open problem.

Our final goal is to study the connection between this Fueter-Sce-Qian CST and $U_a$ and $U_s$ respectively, thus fully answering the question raised in \cite{MR3556034}. As we will show in Section \ref{CSTSect}, the plane wave decomposition of the Fueter-Sce-Qian mapping, provided in Theorem \ref{FullRelRadThm}, allows for a complete description of the connections among all of these mappings. These results are summarized in Theorem \ref{CSTDiagComThm}  and diagram (\ref{CSTDiagComThm}).

Summarizing, in this paper we solve the following problems:
\begin{itemize}
\item[{\bf P1}] establish the connection between the Fueter-Sce-Qian and the generalized CK-extension theorems
\item[{\bf P2}] obtain a plane wave decomposition for the Fueter-Sce-Qian extension map
\item[{\bf P3}] construct a monogenic CST through the Fueter-Sce-Qian theorem and establish its connection with the axial and slice CSTs studied in  \cite{MR3556034}.
\end{itemize}

The plan of the paper is as follows. In Section \ref{prem}, we provide a brief introduction to Clifford analysis. We pay particular attention to the notions of slice and axial monogenic functions and how they can be constructed from real-analytic functions in $\R$ by means of the slice and GCK extensions respectively. In Section \ref{SecFuet}, we review some of the standard facts of the Fueter-Sce-Qian theorem as a {\it bridge} connecting the modules of slice and axial monogenic functions. Section \ref{FuetGCKSect} is devoted to  establishing the relation between the  Fueter-Sce-Qian and the GCK extension theorem, solving thus problem {\bf P1}. In Section \ref{RadTSect}, we study yet another bridge between slice and axial monogenic functions: the dual Radon transform. In particular, we show how to use this transform to decompose the Fueter-Sce-Qian extension map into plane waves, providing a solution to {\bf P2}. Finally, in Section \ref{CSTSect}, we use the previous results to construct a CST using the Fueter-Sce-Qian theorem and study its connection with the other CSTs already defined in the monogenic setting, which solves {\bf P3}.

\section{Preliminary results}\label{prem}
In this section we fix some notations and briefly recall some definitions and results from Clifford analysis. Let  $\mathbb{R}_m$ denote the real associative Clifford algebra with generators $e_1,...,e_m$ satisfying the defining relations $e_je_{\el}+e_{\el}e_j=-2\delta_{j\el}$, for $j,\el=1,...,m$, where $\delta_{j\el}$ is the well-known Kronecker symbol. Every element $a\in \mathbb{R}_m$ can be written in the form 
 \[a=\sum_{A\inc M}a_Ae_A, \hspace{1cm} a_A\in \R,\]
where 
$M:=\{1,\ldots, m\}$ and for any multi-index $A=\{j_1, \ldots, j_k\}\subseteq M$ with $j_1< \ldots < j_k$ we put $e_A:=e_{j_1}\cdots e_{j_k}$ and $|A|=k$. Every $a\in\mathbb{R}_m$ admits a multivector decomposition 
\[a=\sum_{k=0}^m [a]_k, \hspace{.5cm} \mbox{ where } \hspace{.5cm} [a]_k=\sum_{|A|=k}a_A e_A.\]
Here $[\cdot]_k:\R_m\fd \R_m^{(k)}$ denotes the canonical projection of $\R_m$ onto the space of $k$-vectors $\R_m^{(k)}=\textup{span}_\R\{e_A:|A|=k\}$ . Note that $\R_m^{(0)}=\R$, while the space of $1$-vectors $\R_m^{(1)}$ is isomorphic to $\R^m$. Note also that $\R_1\iso\C$ and  $\R_2\iso\Ha$.

There are several automorphisms on $\R_m$ leaving the above multivector structure invariant. One important example is given by the Clifford conjugation $\ba{\cdot}$ defined by
\[
\ba{ab} = \ba{b}\, \ba{a}; \;\;\; \ba{a_A e_A} = a_A \ba{e_A}; \;\;\; \ba{e_j} = -e_j; \;\;\;\; a,b\in\R_m, \; a_A\in\R.
\]
The complex Clifford algebra with the same generators $e_1,...,e_m$ is denoted by $\C_m$ and it is defined by $\C_m=\C\otimes \R_m = \R_m+i\R_m$, where $i\in\C$ is the complex imaginary unit. All concepts introduced so far for $\R_m$ can be reformulated in the case of $\C_m$. In the case of the conjugation, $\C_m$ allows for the so-called Hermitean conjugation $\cdot^\dagger$ defined by
\begin{equation}\label{HermConj}
(a+ib)^\dagger = \ba{a} - i \ba{b}, \;\;\;\;\;\; a,b\in\R_m.
\end{equation}
An important subspace of the real Clifford algebra $\R_m$ is $\R_m^{(0)}\oplus \R_m^{(1)}$ whose elements are called paravectors. This space can be naturally identified with $\R^{m+1}$ and we will denote its elements by
\[
x = x_0 + \m{x} = x_0 + \sum_{j=1}^m x_j e_j, 
\]
where $x_0\in\R$ and $\m{x}:=\sum_{j=1}^m x_j e_j \in\R^{m}$.
It is easily seen that $x \ba{x} = ( x_0 + \m{x}) (x_0 - \m{x})=|x|^2$ where $|x|:=\left(\sum_{j=0}^m x_j^2\right)^{\frac{1}{2}}$ is the Euclidean norm in $\R^{m+1}$. Similarly, for any pair of vectors $\m{x}, \m{y}\in \R^m$ one has  that 
\[
\m{x}\m{y} + \m{y}\m{x}  = -2 \sum_{j=1}^m x_j y_j = -2 \langle\m{x}, \m{y}\rangle,
\]
where $\langle\m{x}, \m{y}\rangle$ denotes the Euclidean inner product in $\R^m$, and therefore, $\m{x}^2=-|\m{x}|^2$ where $|\m{x}|^2 := \sum_{j=1}^m x_j^2$ is the square of the Euclidean norm of $\m{x}$ in  $\R^m$.

Throughout this paper we will mainly consider real-analytic functions with values on $\R_m$ or in its complexification $\C_m$. More precisely, we will work with modules of the form $\mathcal A(E) \otimes \R_m$ or $\mathcal A(E) \otimes \C_m$, where $E\inc\R^\el$ is an open set and $\mathcal A(E)$ denotes the space of $\R$-valued real analytic functions in $E$. Within these modules, we shall focus our attention on two well-known extensions of the notion of holomorphic functions to higher dimensions, the so-called {\it slice monogenic functions} and {\it axially monogenic functions}.

\subsection{Slice monogenic functions}\label{SMFSect}\phantom{a}\\
We now recall the definition and some important properties of slice monogenic functions, see e.g. \cite{MR2520116, MR2742644, MR2684426}. Let us denote by $\Sa^{m-1}$ the $(m-1)$-dimensional sphere of unit $1$-vectors in $\R^{m+1}$, i.e.
\[\Sa^{m-1} :=\left\{\m{x}= \sum_{j=1}^m x_j e_j: |\m{x}| =1\right\}.\]
The vector space $\C_{\m{w}}=\{u+v\, \m{w}: u,v\in \R\}$, passing through 1 and $\m{w}\in \Sa^{m-1}$, is a $2$-dimensional real subspace of $\R^{m+1}$ isomorphic to the complex plane. The isomorphism turns out to be an algebra isomorphism.

\begin{defi}[Slice monogenic functions]
Let $\Om\inc \R^{m+1}$ be an open set. A function $f: \Om \fd \R_m$ is said to be (left) slice monogenic if, for any $\m{w}\in \Sa^{m-1}$, the restrictions of $f$ to the complex planes $\C_{\m{w}}$ are holomorphic on $\Om\cap \C_{\m{w}}$, i.e.\
\[(\pa_u + \m{w}\,  \pa_v) f(u+v\, \m{w})=0.\]
We denote by $\mathcal{SM}(\Om)$ the right $\R_m$-module of slice monogenic functions on $\Om$.
\end{defi}

Slice monogenic functions have good properties when defined on so-called {\it axially symmetric slice domains}.

\begin{defi}[Axially symmetric slice domain]
An open set $\Om$ of $\R^{m+1}\iso \R\oplus \R^m$ is said to be an {\it axially symmetric slice domain} if it satisfies the following conditions:
\begin{itemize}
\item[$i)$]  $\Om\cap \R \neq \emptyset$,
\item[$ii)$] $\Om\cap \C_{\m{w}}$ is a connected open set (i.e.\ a domain) in $\C_{\m{w}}$ for all $\m{w}\in \Sa^{m-1}$,
\item[$iii)$] $\Om$ is SO$(m)$-invariant with respect to the $\R$-axis, i.e.\ $x_0+\m{x} \in \Om \implies x_0+M \m{x} \in \Om$ for all $M\in \textup{SO}(m)$
\end{itemize}
If $m=1$, this definition reduces to connected open sets in the complex plane which are symmetric with respect to the real-axis. In this case we say that $\Om$ is an {\it intrinsic complex domain}. 

An alternative definition for axially symmetric slice domains can be given in terms of intrinsic complex domains. Indeed, $\Om\inc\R^{m+1}$ is an axially symmetric slice domain if and only if there exist an intrinsic domain $\Om_2\inc\C$ such that
\[
\Om=\Om_2\times \Sa^{m-1} =\{u+v\m{\omega}: (u,v) \in \Om_2, \, \m{\omega}\in \Sa^{m-1} \}.
\]
\end{defi}
\begin{rem}
An open set satisfying only conditions $i)$ and $ii)$ is usually called a slice domain, while an open set satisfying only condition $iii)$ is called an axially symmetric open set.
\end{rem}

In the sequel, we will consider slice monogenic functions defined only on axially symmetric slice domains. This is due to the good properties these functions possess in such domains as stated in the following representation theorem. 
\begin{thm}[Representation formula \cite{MR2742644}]
Let $\Om\inc\R^{m+1}$ be an axially symmetric slice domain, and let $\Om_2\inc\C$ be the intrinsic complex domain such that $\Om=\Om_2\times \Sa^{m-1}$. Then any function $f\in \mathcal{SM}(\Om)$ can be written as 
\[
f(u+v\, \m{w}) = \al(u,v)+ \m{w}\, \be (u,v),
\]
where $\al,\be$ are $\R_m$-valued real analytic functions defined in $\Om_2$, such that 
\begin{equation}\label{cond1}
\al(u,v)= \al(u,-v),  \;\;\;\;    \be(u,v) = -\be(u,-v) \;\;\;\;\; \mbox{for all} \;\; (u,v)\in\Om_2,
\end{equation}
and moreover satisfy the Cauchy-Riemman system
\[
\begin{cases}
\pa_u \al - \pa_v \be =0, \\
\pa_v \al +  \pa_u \be =0.
\end{cases}
\] 
\end{thm}

\begin{paragraph}{\bf Extension from holomorphic to slice monogenic.}
The above result shows how to induce slice monogenic functions from so-called {\it holomorphic intrinsic functions}.
\begin{defi}[Holomorphic intrinsic function]\label{IntHolDef}
A holomorphic function $f(z) = \al(u,v)+ i \be (u,v)$ is said to be {\it intrinsic} if it is defined in an intrinsic complex domain $\Om_2$ and satisfies that $f(z)^c = f(z^c)$, where $\cdot^c$ denotes the complex conjugation. This condition means that the real and imaginary part of $f$ satisfy the condition (\ref{cond1}), or equivalently, that the restriction $f|_\R$ of $f$ to the real line is $\R$-valued. 

\noindent We denote by $Hol(\Om_2)$ the space of holomorphic complex functions on $\Om_2$ and by $\mathcal{H}(\Om_2)$ the (real) vector subspace of $Hol(\Om_2)$ complex holomorphic intrinsic functions. In other words,
\begin{align*}
\mathcal{H}(\Om_2) &=\{f\in Hol(\Om_2):  \al(u,v)= \al(u,-v),  \; \be(u,v) = -\be(u,-v)\} \\
&= \{f\in Hol(\Om_2): f|_\R \mbox{ is } \R\mbox{-valued}\}.
\end{align*}
\end{defi}

Clearly, intrinsic holomorphic functions in $\Om_2\in\C$ induce slice monogenic functions on $\Om=\Om_2\times \Sa^{m-1}$ by means of the extension map:
\begin{align*}
S_\C&: \mathcal{H}(\Om_2) \otimes \R_m \fd  \mathcal{SM}(\Om),  &  \al(u,v)+ i \be (u,v) & \; \mapsto  \;\al(x_0,|\m{x}|)+ \m{\omega} \be (x_0,|\m{x}|),
\end{align*}
which consists of replacing the complex variable $z=u+iv$ by the paravector variable $x=x_0+\m{x}$, where the complex unit $i$ is replaced by the unit vector $\m{\omega}=\frac{\m{x}}{|\m{x}|}$ in $\R^m$.

The resulting slice monogenic functions $f(x_0+\m{x})=\al(x_0,|\m{x}|)+ \m{\omega} \be (x_0,|\m{x}|)$ is well-defined at $\m{x}=0$ (independently of $\m{w}$) since $\be(x_0, |\m{x}|)$ is an odd function in the second variable $v=|\m{x}|$. Indeed, using the Taylor expansion of $\be$ we obtain, 
\[
\be(x_0, |\m{x}|) = \sum_{j=0}^\infty \frac{|\m{x}|^{2j+1}}{(2j+1)!} \, \pa_{|\m{x}|}^{2j+1}[\be](x_0,0),
\]
which implies that $\displaystyle \frac{\m{x}}{|\m{x}|} \, \be(x_0, |\m{x}|) \Big|_{\m{x}=0}  = \sum_{j=0}^\infty \frac{ \m{x} |\m{x}|^{2j}}{(2j+1)!} \Big|_{\m{x}=0} \, \pa_{|\m{x}|}^{2j+1}[\be](x_0,0)=0$.
\end{paragraph}

\begin{paragraph}{\bf Extension from real-analytic to slice monogenic.}
The above extension can be realized in terms of real-analytic functions on the real line. This is due to the fact that holomorphic intrinsic functions are uniquely determined by their ($\R$-valued) restrictions to $\R$. Indeed, using the Taylor expansion of $f\in \mathcal{H}(\Om_2)$ around a point $(u,0)$ on the real line, we have
\begin{equation}\label{HolExt}
f(u+iv) = \sum_{j=0}^\infty \frac{(iv)^j}{j!} \, f^{(j)}(u) = \sum_{j=0}^\infty \frac{(iv)^j}{j!} \, \pa_u^{j}[\al] (u,0),
\end{equation}
where $f(u+iv)  = \al(u,v)+ i \be (u,v)$. It is thus clear that $f$ is the unique holomorphic extension of the function $f_0(u)=\al (u,0)$. Moreover, the real and imaginary parts of $f$ are given by 
\begin{align}\label{HolExt1}
\al(u,v)&= \sum_{j=0}^\infty \frac{(-1)^j v^{2j}}{(2j)!} \, \pa^{2j}_u [\al] (u,0), && \mbox{ and } & \be(u,v)&= \sum_{j=0}^\infty \frac{(-1)^j v^{2j+1}}{(2j+1)!} \, \pa^{2j+1}_u [\al] (u,0).
\end{align}
Let us then consider the subset of the real line $\Om_1:=\Om_2\cap \R$ and (abusing the notation) let us denote by $\mathcal{A}(\Om_1)$ the space of $\R$-valued real-analytic functions on $\Om_1$ with unique holomorphic extensions to $\Om_2$. We have thus arrived to the holomorphic extension map $C = \exp(iv\pa_u)$, i.e.\
\begin{align}\label{holExt}
C&: \mathcal{A}(\Om_1) \fd \mathcal{H}(\Om_2), & f_0(u) &\mapsto \sum_{j=0}^\infty \frac{(iv)^j}{j!} f_0^{(j)}(u),
\end{align}
which in turn yields the slice monogenic extension map $S=S_\C \circ C = \exp(\m{x}\pa_{x_0})$, i.e.\
\begin{align}\label{sliceExt}
S&: \mathcal{A}(\Om_1)\otimes \R_m \fd  \mathcal{SM}(\Om), & f_0(x_0) &\mapsto \sum_{j=0}^\infty \frac{\m{x}^j}{j!} f_0^{(j)}(x_0),
\end{align}
\end{paragraph}
These results are all summarized in the following theorem

\begin{thm}[Slice extension theorem]\label{SME}
Under the conditions stated above, the following right $\R_m$-modules are isomorphic,
\[
\mathcal{SM}(\Om) \; \iso \; \mathcal{A}(\Om_1) \otimes \R_m \; \iso \; \mathcal{H}(\Om_2) \otimes \R_m.
\]
In fact, the following diagram commutes
\[
\begin{tikzcd}[row sep = 2em, column sep = 5em]
\mathcal{A}(\Om_1)\otimes \R_m \arrow[r, "C", ] \arrow[dr, "S", labels=below] & \mathcal{H}(\Om_2) \otimes \R_m \arrow[d, "S_\C"] \\
               & \mathcal{SM}(\Om)
\end{tikzcd}
\]
where the map $S=S_\C\circ C =\exp(\m{x}\pa_{x_0})$, given by $S[f_0](x) = \sum_{j=0}^\infty \frac{\m{x}^j}{j!} \, f_0^{(j)}(x_0)$,
is inverted by the restriction operator to the real line, i.e.\ $S[f_0]\big|_{\m{x}=0}=f_0$.
\end{thm}

\subsection{Axial monogenic functions}\phantom{a}\\
Another important extension of the holomorphic function theory to higher dimensions is provided by the notion of monogenic functions. These are functions in the kernel of the  Cauchy-Riemann operator in $\R^{m+1}$ \[
\mathcal{D}_x = \pa_{x_0}+\partial_{\m{x}} = \pa_{x_0}+\sum_{j=1}^m e_j\partial_{x_j}, \]
where the vector (or gradient) part $\partial_{\m{x}} := \sum_{j=1}^m e_j\partial_{x_j}$ is the well-known Dirac operator in $\R^m$.
\begin{defi}[Monogenic functions]\label{MonDef}
A continuously differentiable $\R_m$-valued function $f$, defined on an open set $\Om\inc \R^{m+1}$, is called (left) monogenic if $\mathcal{D}_x f  = 0$ on $\Om$.
\end{defi}
It is easily seen that for $m=1$, monogenic functions on $\R^2$ correspond to holomorphic functions of the variable $x_0+e_1x_1$. Monogenicity also constitutes a refinement of harmonicity since the Cauchy-Riemann operator $\mathcal{D}_x$ factorizes the Laplace operator in $\R^{m+1}$, i.e.\
\begin{align*}
 \Del & := \sum_{j=0}^m \partial_{x_j}^2= \mathcal{D}_x \ba{\mathcal{D}}_x = \ba{\mathcal{D}}_x\mathcal{D}_x.
\end{align*}
Therefore the right $\R_m$-module of monogenic functions defined in an open set $\Om\inc\R^{m+1}$ clearly is a submodule of $\mathcal A(\Om) \otimes \R_m$. Standard references on this setting are \cite{MR1130821, MR697564, MR1169463}. 

Slice monogenic functions can be linked with monogenic functions through the notion of axially monogenic function.

\begin{defi}[Axial monogenic functions]\label{AxMonDef}
Let $\Om\inc\R^{m+1}$ be an axially symmetric slice domain, and let $\Om_2\inc\C$ be the intrinsic complex domain such that $\Om=\Om_2\times \Sa^{m-1}$. A function $f\in\mathcal{A}(\Om)\otimes \R_m$ is said to be axial monogenic if it is monogenic, i.e. $\mathcal{D}_x f =0$, and it has the form
\begin{equation}\label{AxForm}
 f(x_0+\m{x}) = A(x_0, |\m{x}|)+ \m{w} \, B(x_0, |\m{x}|), \;\;\;\;\;\;\; \m{w}=\frac{\m{x}}{|\m{x}|},
\end{equation}
where $A,B\in \mathcal{A}(\Om_2)\otimes \R_m$ satisfy the conditions (\ref{cond1}). We denote by $\mathcal{AM}(\Om)$ the right $\R_m$-module of axially monogenic functions on $\Om$.
\end{defi}
\begin{rem}
As discussed earlier, the fact that $A,B$ satisfy (\ref{cond1}) means that $ f(x_0+\m{x})$ is well-defined at $\m{x}=0$. As a matter of fact, using the Taylor expansions of $A$ and $B$ around points in the real line, we can write $f$ as a power series of the vector variable $\m{x}$, i.e.\
\[
 f(x_0+\m{x}) = \sum_{j=0}^\infty \m{x}^j f_j(x_0), 
\]
where 
$\displaystyle A(x_0, |\m{x}|) = \sum_{j=0}^\infty |\m{x}|^{2j} (-1)^j f_{2j}(x_0)$ and $\displaystyle B(x_0, |\m{x}|) = \sum_{j=0}^\infty |\m{x}|^{2j+1} (-1)^j f_{2j+1}(x_0)$.
\end{rem}

As in the case of slice monogenic functions, there exist an isomorphism between the modules of real-analytic $\R_m$-valued functions on the real line and the module of axially monogenic functions in $\R^{m+1}$. This result, is a particular case of the generalized CK-extension, see \cite{MR1169463}.

\begin{thm}[Generalized CK-extension] \label{GCKthm}
Let $ \Omega_1 \subset \mathbb{R}$ be a real domain, and consider an analytic function $f_0(x_0)\in\mathcal{A}(\Om_1)\otimes\R_m$. Then there exists a unique sequence $ \{f_j(x_0)\}_{j=1}^\infty \inc \mathcal{A}(\Om_1)\otimes\R_m$  such that the series
\[ f(x_0, \underline{x})= \sum_{j=0}^\infty \underline{x}^j f_j(x_0) \]
converges in an axially symmetric slice $(m+1)$-diemensional neighborhood $\Om\inc\R^{m+1}$ of $\Om_1$ and its sum is monogenic, i.e. $ (\partial_{x_0}+ \partial_{\underline{x}}) f(x_0, \underline{x})=0$. Moreover,
\begin{equation}\label{GCK}
 f(x_0, \underline{x})= \Gamma \left(\frac{m}{2}\right) \left(\frac{|\underline{x}| \partial_{x_0}}{2}\right)^{- \frac{m}{2}} \left[ \frac{|\underline{x}| \partial_{x_0}}{2}J_{\frac{m}{2}-1}(|\underline{x}|\partial_{x_0})+ \frac{\underline{x} \partial_{x_0}}{2} J_{\frac{m}{2}}(|\underline{x}|\partial_{x_0})\right] f_0(x_0),
 \end{equation}
where $J_\nu$ is the Bessel function of the first kind of order $\nu$. 

The function in equation (\ref{GCK}) is known as the Genneralized CK-extension of $f_0$, and it is denoted it by  GCK$[f_0](x_0, \underline{x})$. This extension operator defines an isomoprhism between right modules:
 \[\textup{GCK}: \mathcal{A}(\Om_1)\otimes\R_m \fd \mathcal{AM}(\Om),\] 
whose inverse is given by the restriction operator to the real line, i.e.\ GCK$[f_0](x_0,0)=f_0(x_0)$.
\end{thm}

\section{Fueter-Sce-Qian theorem: a bridge}\label{SecFuet}
We now turn our attention to the connection between the modules of slice monogenic and axial monogenic functions. This connection is provided by the well-known Fueter-Sce-Qian theorem. 

Initially, this relation was established in the quaternionic case by Fueter \cite{MR1509515}, and later generalized by Sce \cite{MR97386} to the Euclidean space $\R^{m+1}$ for odd values of the dimension $m$. Specifically, this initial result showed that the pointwise differential operator
\[
\tau_m = \Del^{\frac{m-1}{2}} \;\;\;\;\;\; (\mbox{where } m\in\N \mbox{ is odd})
\]
maps slice monogenic functions into axial monogenic ones.

Later in \cite{MR1485323}, Qian extended this result to any dimension $m\in\N$ with the help of the Fourier multiplier 
\begin{equation}\label{FourMult}
\left(-\Del\right)^{\frac{m-1}{2}}= \mathcal{F}^{-1} \left(2\pi|\cdot|\right)^{m-1} \mathcal{F},
\end{equation} 
which gives meaning to the fractional powers of the Laplacian corresponding to the cases where $m$ is even. Here $\mathcal{F}$ and $\mathcal{F}^{-1}$ are the Fourier and inverse Fourier transform in $\R^{m+1}$ given, respectively, by
\begin{align*}
\mathcal{F}[\phi](\xi) &= \int_{\R^{m+1}} e^{2\pi i \langle x,\xi\rangle} \phi(x)\, dx && \mbox{ and }& \mathcal{F}^{-1}[\phi](\xi) &= \int_{\R^{m+1}} e^{-2\pi i \langle x,\xi\rangle} \phi(x)\, dx.
\end{align*}
Qian's extension makes use of the following more constructive approach. Consider a meromorphic intrinsic function $f(z)$ on $\C$, together with its Laurent expansion around $z=0$
\[
f(z) = \sum_{j\in\Z} a_j z^j.
\]
It is clearly seen that all the coefficients $a_j$ in the above expansion must be real numbers and therefore, they play no role when applying the slice extension map $S_\C$. Thus, to determine the action of the Fueter-Sce-Qian mapping $\tau_m$ on the slice extension $f(x_0+\m{x}):=S[f](x)$ of $f$, it is enough to consider only the actions $\tau_m\left[x^j\right], \, j\in\Z$. 

A direct computation of these actions, using the Fourier multiplier definition (\ref{FourMult}), is possible for all dimensions $m\in \N$ only when we consider negative powers of $x$, i.e.\ $x^{-\el}$ with $\el\in\N$. Observe that $x^\el$ ($\el\in\N$) is not in the Schwartz class of rapidly decreasing functions and therefore,  $\Del^{\frac{m-1}{2}}\left[x^\el\right]$ is not well-defined as a function if $m$ is even. However, it is possible to define suitable actions of the Fueter-Sce-Qian mapping on positive powers of the paravector $x$ by means of the Kelvin inversion 
\begin{equation}\label{KIDef}
I[f](x) =\frac{\ba{{x}}}{|x|^{m+1}} f\left(\frac{\ba{x}}{|x|^{2}}\right),
\end{equation}
which maps monogenic functions into monogenic functions.

In general, the following relations have been established for the action of the Fueter-Sce-Qian mapping on integer powers of $x$ for all dimensions $m\in\N$,
\begin{equation}\label{FQTMap}
\tau_m\left[x^\el\right] = \begin{cases}
\sgn(-x_0)^{m-1} \Del^{\frac{m-1}{2}}\left[x^\el\right], & \el<0, \\
0, & 0\leq \el \leq m-2, \\
I\left[\Del^{\frac{m-1}{2}}\left[x^{-\el+m-2}\right]\right], & m-1\leq \el.
\end{cases}
\end{equation}
In \cite{MR1485323}, it was shown that (\ref{FQTMap}) indeed provides a suitable extension of the pointwise differential operator $\Del^{\frac{m-1}{2}}\left[x^\el\right]$ when $m$ is odd since, in this case, both approaches coincide. 
\begin{rem}
It is worth noticing that the above definition of the Fueter-Sce-Qian mapping is slightly different form the customary definition used in the literature, see e.g.\ \cite{MR1485323, MR1957490}. Our modification consists in the introduction of the (almost constant) factor  $\sgn(-x_0)^{m-1}$ after the action  $\Del^{\frac{m-1}{2}}\left[x^\el\right]$ when $\el<0$. As we will show in the following section, this small difference allows for a better description of $\tau_m$ in terms of the GCK operator. Observe also that this change still preserves the property that $\tau_m$ extends the pointwise differential operator $\Del^{\frac{m-1}{2}}$ from $m$ odd to $m\in\N$. Indeed, it is obvious that $\sgn(-x_0)^{m-1}\equiv 1$ if $m$ is odd.
\end{rem}

More explicit expressions for the actions  (\ref{FQTMap}) will be provided in section \ref{EvenCaseSection}, see also \cite{MR1485323, MR1957490, MR3539492}. We can now summarize the Fueter-Sce-Qian theorem as follows.

\begin{thm}[Fueter-Sce-Quian theorem]\label{FueterThm1}
Let $f(u+iv)=\al(u,v)+i\be(u,v)$ be an intrinsic holomorphic function defined on an intrinsic complex domain $\Om_2\inc \C$, and put $\m{w}=\frac{\m{x}}{r}$ with $r=|\m{x}|$. Then,
\begin{equation}\label{FuePrim}
\tau_m \left[f(x_0+\m{x})\right] = \tau_m \left[ \al(x_0,r)+\m{w}\, \be(x_0,r)\right],
\end{equation}
is axially monogenic in the axially symmetric slice domain $\Om= \Om_2\times \Sa^{m-1} = \{(x_0,\m{x})\in \R^{m+1}: (x_0,|\m{x}|)\in \Om_2\}$. 


\end{thm}
Since every intrinsic holomorphic function $f$ is the unique holomorphic extension of a real-analytic function $f_0$ on the real line, Fueter's theorem can be rewritten in the following terms.
\begin{thm}[Fueter-Sce-Quian theorem] \label{FSQThm}
Let $\Om_1\inc \R$ be a real domain and $f_0\in\mathcal{A}(\Om_1)\otimes \R_m$. Then $\tau_m \circ S [f_0](x_0,\m{x})$ is an axial monogenic function on a $(m+1)$-dimensional axially symmetric slice neighborhood $\Om\inc \R^{m+1}$ of  $\Om_1$.
\end{thm}

The Fueter-Sce-Qian theorem establishes a link between slice monogenic functions and axially monogenic functions. Indeed, combining Theorems \ref{SME} and \ref{FSQThm} yields the mapping property \[\tau_m: \mathcal{SM}(\Om)\fd \mathcal{AM}(\Om).\]
The actions of the Fueter-Sce-Qian, the GCK and the slice monogenic extension maps can be summarized in the following diagram:

\begin{equation}\label{Diag1}
\begin{tikzcd}[row sep = 3em, column sep = 6em]
\mathcal{A}(\Om_1)\otimes \R_m \arrow[r, "S", rightarrow] \arrow[d, dashed, "?"] & \mathcal{SM}(\Om) \arrow[d, "\tau_m"] \\
\mathcal{A}(\Om_1)\otimes \R_m \arrow[r, "\textup{GCK}", rightarrow]                 & \mathcal{AM}(\Om)
\end{tikzcd}
\end{equation}
Obviously, the map GCK does not coincide with the Fueter-Sce-Qian map $\tau_m\circ S$ since GCK is an isomorphism between right modules while $\tau_m\circ S$ is not. Nevertheless, the above diagram can be completed (and made commutative) by adding the missing left vertical arrow. We address this problem in the next section.

\section{Fueter-Sce-Qian's theorem and generalized CK-extension}\label{FuetGCKSect}
In this section, we uncover the relation between the Fueter map and the generalized CK-extension, hence making the diagram (\ref{Diag1}) complete and commutative. First we study the case where the dimension $m\in\N$ is odd, and $\tau_m=\Del^{\frac{m-1}{2}}$ is a pointwise differential operator. Later we focus on the case where $m$ is even and $\tau_m$ is defined with the help of the Fourier multiplier and the Kelvin inversion as in (\ref{FQTMap}). 

\subsection{The odd dimensional case}\phantom{a}\\
In this case, the action of the pointwise differential operator $\Del^{\frac{m-1}{2}}$ on slice monogenic functions has been explicitly computed, see e.g. \cite[Lem.3.2]{DixanThesis} and \cite[Thm. 11.33]{MR2369875}. This result will prove extremely useful when deriving the connection between the Fueter-Sce-Qian's theorem and the generalized CK-extension.

\begin{lem}[\cite{DixanThesis, MR2369875}]\label{ExpFuetExpr}
If $m\in\N$ is odd, then the action of the pointwise differential operator $\tau_m=\Del^{\frac{m-1}{2}}$ on a slice monogenic function
\[
f(x_0+\m{x})= \al(x_0,r)+\m{w}\, \be(x_0,r),  \;\;\;\;\;\;\; \mbox{ with } \;\; r=|\m{x}| \; \mbox{ and } \;\; \m{w}=\frac{\m{x}}{r},  
\]
is given by 
\[
\tau_m \left[f(x_0+\m{x})\right] = A(x_0,r) + \m{w}\, B(x_0,r),
\]
with 
\begin{align}\label{ExplAB}
A(x_0,r) &= (m-1)!! \left(\frac{1}{r}\pa_r\right)^{\frac{m-1}{2}}[\al](x_0,r), & &\mbox{and} & B(x_0,r) &= (m-1)!! \left(\pa_r \frac{1}{r}\right)^{\frac{m-1}{2}}[\be](x_0,r).
\end{align}
\end{lem}
 We can now formulate the main result of this section.
\begin{thm}\label{res1}
Let $f(u+iv)=\al(u,v)+i\be(u,v)$ be an intrinsic holomorphic function defined on an intrinsic complex domain $\Om_2\inc \C$. Then for $m$ odd and $r=|\m{x}|$we have
\begin{eqnarray*}
\Delta^{\frac{m-1}{2}}\left[f(x_0+ \underline{x})\right]&=&(m-1)!! \; \textup{GCK} \left[\left( \frac{1}{r} \partial_r\right)^{\frac{m-1}{2}}[\al](x_0,0) \right]\\
&=& (-1)^{\frac{m-1}{2}} \frac{(m-1)!!}{(m-2)!!} \; \textup{GCK}\left[f^{(m-1)}(x_{0})\right].
\end{eqnarray*}
Setting $\gamma_m= (-1)^{\frac{m-1}{2}} \frac{(m-1)!!}{(m-2)!!}$ and $\Om_1=\Om_2\cap \R$, we obtain the following commutative diagram.
\begin{equation}\label{Diag2}
\begin{tikzcd}[row sep = 3em, column sep = 6em]
\mathcal{A}(\Om_1)\otimes \R_m \arrow[r, "S", rightarrow] \arrow[d, "\gamma_m \pa_{x_0}^{m-1}", labels=left] & \mathcal{SM}(\Om) \arrow[d, "\Del^{\frac{m-1}{2}}"] \\
\mathcal{A}(\Om_1)\otimes \R_m \arrow[r, "\textup{GCK}", rightarrow]                 & \mathcal{AM}(\Om)
\end{tikzcd}
\end{equation}
\end{thm}
\begin{proof}
The fact that $f$ is an holomorphic intrinsic function means that $\al(u,v)$ is an even analytic function in $v$ while $\be(u,v)$ is an odd analytic function in $v$. Thus, after setting $u=x_0$ and $v=r$, we easily see that the functions
\begin{align*}
A(x_0,r) &= (m-1)!! \left(\frac{1}{r}\pa_r\right)^{\frac{m-1}{2}}[\al](x_0,r), & &\mbox{and} & B(x_0,r) &= (m-1)!! \left(\pa_r \frac{1}{r}\right)^{\frac{m-1}{2}}[\be](x_0,r),
\end{align*}
also are even and odd analytic functions in the variable $r$ respectively.  Using Lemma \ref{ExpFuetExpr} we identify these functions as the components of the axial monogenic function that results from the action of the Fueter-Sce-Qian map on $f$, i.e.\
\[
\Del^{\frac{m-1}{2}} f(x_0+\m{x}) = A(x_0,r)+\m{w} B(x_0, r). 
\]
By virtue of the generalized CK-extension Theorem \ref{GCKthm}, the above axially monogenic function is completely determined by its restriction to the real line 
\[\Delta^{\frac{m-1}{2}} f(x_0+ \underline{x})\big|_{\underline{x}=0}.\] 
Since $B(x_0, r)$ is odd in the variable $r$, we easily obtain that 
\[
\Delta^{\frac{m-1}{2}} f(x_0+ \underline{x})\Big|_{\underline{x}=0} = A(x_0, 0) = (m-1)!! \left(\frac{1}{r}\pa_r\right)^{\frac{m-1}{2}}[\al](x_0,0).
\]
Therefore, 
\[
\Delta^{\frac{m-1}{2}} f(x_0+ \underline{x}) = (m-1)!! \;\textup{GCK}\left[\left(\frac{1}{r}\pa_r\right)^{\frac{m-1}{2}}[\al](x_0,0)\right].
\]
From formula (\ref{HolExt1}) we now recall that 
\[
\al(x_0,r) = \sum_{j=0}^\infty \frac{(-1)^j r^{2j}}{(2j)!} \, \pa^{2j}_{x_0} [\al] (x_0,0), 
\]
which for any $\el\in\N$ yields
\begin{align*}
\left(\frac{1}{r}\pa_r\right)^{\el} [\al](x_0,r) &= \sum_{j=\el}^\infty (-1)^j \frac{(2j)(2j-2)\cdots(2j-2\el+2)}{(2j)!} \, r^{2j-2\el} \, \pa^{2j}_{x_0} [\al] (x_0,0) \\
&= \sum_{j=0}^\infty (-1)^{j+ \ell} 2^\ell \frac{(j+ \ell)!}{j! (2j+2 \ell)!} \,  r^{2 j} \partial_{x_0}^{2j+2 \ell}[\al](x_0,0).
\end{align*}
Taking $r=0$, we obtain
\begin{align*}
\left(\frac{1}{r}\pa_r\right)^{\el} [\al](x_0,0) &=   \frac{(-1)^{ \ell}}{(2 \ell-1)!!} \, \partial_{x_0}^{2 \ell}[\al](x_0,0) =  \frac{(-1)^{ \ell}}{(2 \ell-1)!!} \, f^{(2\el)}(x_0),
\end{align*}
which proves the desired result when substituting $\el=\frac{m-1}{2}$. 
\end{proof}

\subsection{The even case} \label{EvenCaseSection}{\phantom{a}}\newline
We now extend  Theorem \ref{res1} to any dimension $m\in\N$ regardless of the parity of $m$. To that end, we make use of the more constructive approach developed by Qian and his collaborators, see e.g.\ \cite{MR1485323, MR1957490}. 
 As we discussed in section \ref{SecFuet}, this approach focusses on the basic actions $\tau_m\left[x^j\right]$ defined in (\ref{FQTMap}) for all $j\in\Z$. 
These actions on integer powers of the paravector $x$  have been explicitly computed for all dimensions $m\in\N$ and all powers $j\in\Z$, in terms of the so-called monogenic monomials, see e.g.\ \cite{MR1485323, MR1957490, MR3539492}. 
\begin{defi}[Monogenic monomials]\label{DefMM}
 Let $m,k\in\N$, we define the monogenic monomials $P^{(-k)}$ and $P^{(k-1)}$ respectively by 
 \begin{align*}
 P^{(-k)} &:= \frac{(-1)^{k-1}\sigma_{m+1} \lan_m}{(k-1)!} \;  \pa_{x_0}^{k-1} [E], & P^{(k-1)} &:= I[ P^{(-k)} ].
 \end{align*}
Here $\lan_m= 2^{m-1} \, \left(\Gamma\left(\frac{m+1}{2}\right)\right)^2$, $\sigma_{m+1} =\frac{2\pi^{\frac{m+1}{2}}}{\Gamma\left(\frac{m+1}{2}\right)}$ is the surface area of the unit sphere $\Sa^{m}$ in $\R^{m+1}$, $E(x)$ is the so-called Cauchy kernel, i.e.\ the fundamental solution of $\mathcal{D}_x = \pa_{x_0}+\pa_{\m{x}}$ given by
 \[E(x)=\frac{1}{\sigma_{m+1}} \frac{\ba{x}}{|x|^{m+1}},\] 
while $I$ is the Kelvin inversion defined in (\ref{KIDef}), i.e.\
 \[I[f](x)=\sigma_{m+1} \, E(x) \, f\left(\frac{\ba{x}}{|x|^2}\right).\] 
\end{defi}
Using the Fourier multiplier definition of the fractional Laplacian (\ref{FourMult}), it has been shown that (see for example \cite{MR1957490, MR1485323, MR3539492})
\[
\left(-\Del\right)^{\frac{m-1}{2}} [x^{-k}] := \mathcal{F}^{-1} \left[\left(2\pi|\cdot|\right)^{m-1} \mathcal{F}[(\cdot)^{-k}]\right](x) = P^{(-k)}(x)
\]
Combining this identity with the definition of $\tau_m$ given in (\ref{FQTMap}) one easily obtains the following result.
\begin{thm}[{{\cite{MR1485323, MR1957490}}}]
\label{MonThmQ}
Let $m\in\N$, $\el\in\Z$. The actions of the Fueter-Sce-Qian mapping $\tau_m[x^\el]$ defined in (\ref{FQTMap}) are given by the expressions
\[
\tau_m[x^\el] = 
(-1)^{\frac{1-m}{2}}\begin{cases}
\sgn(-x_0)^{m-1} P^{(\el)}(x), & \el<0,\\
0, & 0\leq\el\leq m-2,\\
P^{(\el+1-m)}(x), & m-1\leq \el.
\end{cases}
\]
\end{thm}
From the above theorem, it is clear that our goal now reduces to find how the monogenic monomials can be expressed in terms of the generalized CK-extension map. 
To that end, we first note that the Cauchy kernel $E(x)$ is an axial monogenic function on $\R^{m+1}\setminus\{0\}$ since
\[
E(x_0,\m{x}) = x_0(x_0^2+r^2)^{-\frac{m+1}{2}} - \m{\omega} \, r (x_0^2+r^2)^{-\frac{m+1}{2}},  \;\;\;\;\;\;\; \mbox{ with } \;\; r=|\m{x}| \; \mbox{ and } \;\; \m{w}=\frac{\m{x}}{r}. 
\]
It is thus clear that the monomials with negative indices $P^{(-k)}$, $k\in\N$, are also axial monogenic. The axial monogenicity of the rest of the monomials $P^{(k)}$, with $k\in\N$, can be shown as a consequence of the following result. 
\begin{lem}\label{IvsGCK}
The following statements hold:
\begin{itemize}
\item[$i)$] The Kelvin inversion $I$ preserves axial monogenicity.

\item[$ii)$] Given a domain $\Om_1\inc \R$ and $f_0\in\mathcal{A}(\Om_1) \otimes \R_m$, we have
\begin{equation*}
I \circ \textup{GCK}[f_0] =  \sgn (x_0)^{m+1}\, \textup{GCK} \left[ x_0^{-m} f_0(x_{0}^{-1})\right].
\end{equation*}
\end{itemize}
\end{lem}
\begin{proof}
In \cite[Chpt. II]{MR1169463}, it was proved that the Kelvin inversion $I$ preserves monogenicity. Thus, to prove $i)$, it suffices to show that $I$ preserves the axial form (\ref{AxForm}). 

\noindent We first observe that  $|x_0+\m{x}|^\gam = (x_0^2+ |\m{x}|^2)^\frac{\gam}{2}$ ($\gam\in\R$) is a scalar function of $x_0$ and $|\m{x}|$. Thus for any axial monogenic function $f(x_0,\m{x})= A(x_0, |\m{x}|)+ \m{x} \, B(x_0, |\m{x}|)$ we have that
\begin{align*}
I[f](x_0,\m{x}) &= \frac{x_0-\m{x}}{|x_0+\m{x}|^{m+1}} \left[ A\left(\frac{x_0}{|x_0+\m{x}|^2}, \frac{|\m{x}|}{|x_0+\m{x}|^2}\right)- \frac{\m{x}}{|x_0+\m{x}|^2} \, B\left(\frac{x_0}{|x_0+\m{x}|^2}, \frac{|\m{x}|}{|x_0+\m{x}|^2}\right)\right]\\
&= A_1(x_0, |\m{x}|)+ \m{x} \, B_1(x_0, |\m{x}|),
\end{align*}
for a suitable pair of functions $A_1$ and $B_1$.
\noindent The second statement $ii)$ easily follows from the properties of the generalized CK-extension (see Theorem \ref{GCKthm}). Indeed,
\[
I \circ \textup{GCK}[f_0] (x_0,\m{x}) = \frac{x_0-\m{x}}{|x_0+\m{x}|^{m+1}} \,  \textup{GCK}[f_0] \left(\frac{x_0}{|x_0+\m{x}|^2}, -\frac{\m{x}}{|x_0+\m{x}|^2}\right).
\]
Hence, the restriction of this axial monogenic function to the real line is given by
\begin{align*}
I \circ \textup{GCK}[f_0] (x_0,0) &=  \frac{x_0}{|x_0|^{m+1}} \,  \textup{GCK}[f_0] \left(x_0^{-1}, 0\right) \\
&=   \sgn (x_0)^{m+1} x_0^{-m} f_0(x_{0}^{-1}),
\end{align*}
which proves the result.
\end{proof}

We can now write the monogenic monomials $P^{(k)}$ as generalized CK-extensions of suitable initial analytic functions of one real variable.

\begin{prop}\label{MonMonFGCK}
For all $m,k \in \mathbb{N}$ we have that
\begin{align}
P^{(-k)} &= \frac{\lambda_m(m+k-2)!}{(k-1)! (m-1)!}\, \sgn(x_0)^{m-1} \, \textup{GCK}[x_0^{-k-m+1}], 
\label{resu1}\\
P^{(k-1)} &=  \frac{\lambda_m(m+k-2)!}{(k-1)! (m-1)!}\; \textup{GCK}[ x_0^{k-1}], 
\label{res2}
\end{align}
Or equivalently, 
\begin{align}
P^{(-k)} &= \frac{\lambda_m}{(m-1)!}\, \sgn(-x_0)^{m-1} \, \textup{GCK} \circ \pa_{x_0}^{m-1}[x_0^{-k}], 
\label{res3}\\
P^{(k-1)} &=\frac{\lambda_m}{(m-1)!}\; \textup{GCK} \circ \pa_{x_0}^{m-1} [ x_0^{m+k-2}]. 
\label{res4}
\end{align}
\end{prop}
\begin{proof}
From $ P^{(-k)} = \frac{(-1)^{k-1}\sigma_{m+1} \lan_m}{(k-1)!} \;  \pa_{x_0}^{k-1} [E]$ we obtain 
\begin{align*}
P^{(-k)} \big|_{\m{x}=0} &=  \frac{(-1)^{k-1} \lan_m}{(k-1)!} \; \pa_{x_0}^{k-1}\left[\frac{x_0}{|x_0|^{m+1}}\right] \nonumber \\ 
& =  \frac{(-1)^{k-1} \lan_m}{(k-1)!} \,  \sgn(x_0)^{m+1} \, \pa_{x_0}^{k-1}\left[x_0^{-m}\right] \nonumber  \\
&= \frac{\lan_m (m+k-2)!}{(k-1)! (m-1)!} \,  \sgn(x_0)^{m-1} \, x_0^{-k-m+1}. 
\end{align*}
On the other hand, substituting the indentity $x_0^{-k-m+1} = (-1)^{m-1} \frac{(k-1)!}{(k+m-2)!} \pa_{x_0}^{m-1}[x_0^{-k}]$ in the above equality yields
\[
P^{(-k)} \big|_{\m{x}=0} = \frac{\lambda_m}{(m-1)!}\, \sgn(-x_0)^{m-1} \, \pa_{x_0}^{m-1}[x_0^{-k}].
\]
Using the fact that axially monogenic functions are completely determined by their restrictions to the real line, we obtain from the two last equalities that (\ref{resu1}) and (\ref{res3}) hold.

\noindent For $P^{(k-1)}$ with $k\in\N$ we have, by virtue of Lemma \ref{IvsGCK} and (\ref{resu1}), that 
\begin{align*}
P^{(k-1)} &= I[P^{(-k)}] \\
&= \frac{\lambda_m(m+k-2)!}{(k-1)! (m-1)!}\, \sgn(x_0)^{m-1} \, I\circ \textup{GCK}[x_0^{-k-m+1}] \\
&=  \frac{\lambda_m(m+k-2)!}{(k-1)! (m-1)!}\; \textup{GCK}[ x_0^{k-1}].
\end{align*}
Substituting the identity $x_0^{k-1}=\frac{(k-1)!}{(m+k-2)!} \pa_{x_0}^{m-1}[x_0^{m+k-2}]$ in the above equality we obtain,
\[
P^{(k-1)}  =  \frac{\lambda_m}{(m-1)!}\; \textup{GCK} \circ \pa_{x_0}^{m-1} [ x_0^{m+k-2}],
\]
which completes the proof.
\end{proof}

We can now extend the above relation between the Fueter-Sce-Qian map and the generalized CK-extension when acting on the basic monomials $x_0^k$, $k\in\Z$, to general analytic functions by considering their Laurent expansions.
\begin{thm}\label{GenRel}
Let $ \Omega_2 \subset\C$ be an intrinsic complex domain and let $f: \Omega_2 \to \mathbb{C}$ be a holomorphic intrinsic function. Then for all dimensions $m \in \mathbb{N}$ we have
\begin{equation}
\tau_m \left[f(x_0+ \underline{x}) \right]=  \frac{(-1)^{\frac{1-m}{2}} 2^{m-1} }{(m-1)!}  \Gamma \left( \frac{m+1}{2}\right)^2  \, \textup{GCK} \circ \partial_{x_0}^{m-1}[f(x_{0})].
\end{equation}
Setting $\gamma_m=   \frac{(-1)^{\frac{1-m}{2}} 2^{m-1} }{(m-1)!}  \Gamma \left( \frac{m+1}{2}\right)^2 $ and $\Om_1=\Om_2\cap \R$, we obtain the following extension of the commutative diagram (\ref{Diag2}) to all dimensions $m\in\N$
\begin{equation*}
\begin{tikzcd}[row sep = 3em, column sep = 6em]
\mathcal{A}(\Om_1)\otimes \R_m \arrow[r, "S", rightarrow] \arrow[d, "\gamma_m \pa_{x_0}^{m-1}", labels=left] & \mathcal{SM}(\Om) \arrow[d, "\tau_m"] \\
\mathcal{A}(\Om_1)\otimes \R_m \arrow[r, "\textup{GCK}", rightarrow]                 & \mathcal{AM}(\Om)
\end{tikzcd}
\end{equation*}
\end{thm}

\begin{rem}
The previously defined constant $\gamma_m =   \frac{(-1)^{\frac{1-m}{2}} 2^{m-1} }{(m-1)!}  \Gamma \left( \frac{m+1}{2}\right)^2 $ is an extension to all dimensions $m\in\N$ of the constant $\gamma_m= (-1)^{\frac{m-1}{2}} \frac{(m-1)!!}{(m-2)!!}$ introduced in Theorem \ref{ExplAB} for odd values of $m$. Indeed, if $m$ is odd, then
\[
  \Gamma \left( \frac{m+1}{2}\right)^2 =   \left[\left(\frac{m-1}{2}\right) \left(\frac{m-3}{2}\right) \ldots \left(\frac{2}{2}\right)\right]^2 = \frac{\big((m-1)!!\big)^2}{2^{m-1}},
\]
and $(-1)^{\frac{1-m}{2}}=(-1)^{\frac{m-1}{2}}$ since the power $\frac{m-1}{2}$ is integer. The combination of these two facts easily yields 
\[
 \frac{(-1)^{\frac{1-m}{2}} 2^{m-1} }{(m-1)!}  \Gamma \left( \frac{m+1}{2}\right)^2 = (-1)^{\frac{m-1}{2}} \frac{(m-1)!!}{(m-2)!!}, \;\;\;\mbox{ for } m \mbox{ odd}.
\]
\end{rem}

\begin{proof}
We can assume without loosing generality that $f$ is holomorphic around the origin.  If that is not the case, we can always obtain such a function by applying a translation argument. 
\noindent Let us consider the Laurent expansion of $f$ at $z=0$, i.e.
\[
f(z) = \sum_{j\in\Z} a_j\,  z^j, \;\;\;\; a_j\in\R.
\]
The action of the Fueter mapping on this function is thus given by
\[
\tau_m f(x_0+ \underline{x}) =  \sum_{j\in\Z} a_j \; \tau_m \left[ x^j \right], \;\;\;\; a_j\in\R.
\]
We recall that the action of $\tau_m$ does not affect the convergence of the above series, see \cite{MR1957490}.  
Hence, combining Theorem \ref{MonThmQ} and Proposition \ref{MonMonFGCK}, we obtain
\begin{align*}
\tau_m f(x_0+ \underline{x}) &= \sum_{k=1}^\infty  a_{-k}
\; \tau_m \left[x^{-k} \right]  +  \sum_{k=1}^\infty  a_{k}
\; \tau_m \left[ x^{k} \right] \\
&= (-1)^{\frac{1-m}{2}} \sum_{k=1}^\infty  a_{-k} \, \sgn(-x_0)^{m-1} \;P^{(-k)}  +  (-1)^{\frac{1-m}{2}}\sum_{k=m-1}^\infty  a_{k} \;P^{(k+1-m)}\\
&=(-1)^{\frac{1-m}{2}} \sum_{k=1}^\infty  a_{-k} \, \sgn(-x_0)^{m-1} \;P^{(-k)}  +  (-1)^{\frac{1-m}{2}} \sum_{k=0}^\infty  a_{k+m-1} \;P^{(k)}\\
&= (-1)^{\frac{1-m}{2}}\sum_{k=1}^\infty  a_{-k} \; \frac{\lambda_m}{(m-1)!}\,  \textup{GCK} \circ \pa_{x_0}^{m-1}[x_0^{-k}] \\
&\phantom{=} + (-1)^{\frac{1-m}{2}} \sum_{k=0}^\infty  a_{k+m-1} \;  \frac{\lambda_m}{(m-1)!}\; \textup{GCK} \circ \pa_{x_0}^{m-1} [ x_0^{m+k-1}] \\
&=(-1)^{\frac{1-m}{2}} \frac{\lambda_m}{(m-1)!}\; \textup{GCK} \circ \pa_{x_0}^{m-1}  \left[ \sum_{k=1}^\infty  a_{-k} \; x_0^{-k} + \sum_{k=0}^\infty  a_{k+m-1} \; x_0^{m+k-1} \right]\\
&= (-1)^{\frac{1-m}{2}} \frac{\lambda_m}{(m-1)!}\; \textup{GCK} \circ \pa_{x_0}^{m-1}  \left[  \sum_{k=1}^\infty  a_{-k} \; x_0^{-k} + \sum_{k=0}^\infty  a_{k} \; x_0^{k} \right]\\
&=(-1)^{\frac{1-m}{2}} \frac{\lambda_m}{(m-1)!}\; \textup{GCK} \circ \pa_{x_0}^{m-1}[f(x_0)], 
\end{align*}
which completes the proof.
\end{proof}

\subsection{Connection with Appell sequences of monogenic polynomials}\label{APSect}{\phantom{a}}\newline
We now describe a simple application of Theorem \ref{GenRel} that allows to connect the monogenic monomials $P^{(k)}$ ($k\in\N$) with a well-known Appell sequence of monogenic polynomials. {In the forthcoming paper \cite{F-S_Ant_Kam_Ali}, we will use this connection to study extensions of the classical Fock and Hardy spaces to the Clifford analysis setting.}

In classical terms, a sequence of polynomials $\{P_k\}_{k\in\N}$ of the real variable $x_0$ is called an Appell sequence if it satisfies the property \[\frac{d}{dx_0} P_{k} = k P_{k-1}.\] Important examples of this class of polynomial sequences are given by the classical monomials $\{x_0^k\}_{k\in\N}$ as well as Hermite, Bernoulli and Euler polynomials. 

Appell sequences have also been considered in the Clifford analysis setting with respect to the action on monogenic polynomials of the so-called hypercomplex derivative $\frac{1}{2}(\pa_{x_0}-\pa_{\m{x}})$. The study of Appell sequences of monogenic polynomials has allowed the extension of a wide variety of special functions and function spaces from the holomorphic function theory to the framework of Clifford analysis, see e.g. \cite{MR2769480, MR3994389, MR3057360}.

One important example of an Appell sequence of monogenic polynomials is given by
\begin{equation}\label{MalPol}
Q_k^m(x) = \sum_{j=0}^k T_j^k(m) \,  \ba{x}^{\,k-j} x^j,
\end{equation}
where $T_j^k(m) = \frac{k!}{(m)_k} \frac{\left(\frac{m+1}{2}\right)_{k-j} \left(\frac{m-1}{2}\right)_{j}}{(k-j)! j!}$ and $(a)_ \el$  denotes the Pochammer symbol, i.e. $(a)_ \el=\frac{\Gam(a+\el)}{\Gam(a)}$, $\el\in\N$. The family $\{Q_k^m\}_{k\in\N}$ has been used in the study of Laguerre monogenic polynomials in \cite{MR2769480} and in the study Bargman-Fock type transforms in the quaternionic setting in \cite{MR3994389}. The following properties can be verified:
\begin{itemize}
\item[$i)$] The polynomials  $Q_k^m$ are indeed monogenic, i.e. $(\pa_{x_0}+\pa_{\m{x}}) Q_k^m = 0$.
\item[$ii)$] The polynomials  $Q_k^m$  satisfy the Appell condition $\frac{1}{2}(\pa_{x_0}-\pa_{\m{x}}) Q_k^m  = k Q_{k-1}^m.$
\item[$iii)$] $Q_k^m(1) = 1$ for all $k\in\N$, i.e. 
$\sum_{j=0}^k T_j^k(m) =1$.
\end{itemize}
It is easily seen that each function $Q_k^m(x)=Q_k^m(x_0,\m{x})$ is of axial type (see Definition \ref{AxMonDef}). It is thus clear that  $Q_k^m(x)$ is the generalized CK-extension of its restriction to the real line, i.e.\
\[
Q_k^m(x_0+\m{0}) =  \sum_{j=0}^k x_0^k T_j^k(m) = x_0^k.
\] 
Using the identity $\pa_{x_0}^{m-1}[x_0^{m-1+k}]=\frac{(m+k-1)!}{k!}x_0^{k}$, we arrive at the following corollary of Theorem \ref{GenRel}.
\begin{cor}[Explicit formula for $\tau_m(x^k)$, $k\in\N$]
Let $m,k\in\N$, then
\[\tau_m \left[x^{m-1+k}\right] =  \gamma_m \, \textup{GCK} \circ \pa_{x_0}^{m-1} [x_0^{m-1+k}] = \gamma_m\, \frac{(m-1+k)!}{k!} \,Q_k^m(x).\]
\end{cor}
The above formula was obtained by means of direct computations in the quaternionic case ($m=3$) in \cite{MR3994389}. In the following section we shall provide plane wave expansions for the polynomials $Q_k^m(x)$.

\section{The dual Radon transform: another bridge}\label{RadTSect}
The goal of this section is to uncover the relation between the Fueter-Sce-Qian map $\tau_m$ and the well-known dual Radon transform. The question of whether such a relation exists comes motivated by the fact that the dual Radon transform provides yet another bridge between the modules of slice and axial monogenic functions, see \cite{MR3412341}.  Let us first recall the classical definition of this transform (see \cite{MR754767}) and adapt it to our setting.
\begin{defi}[Dual Radon transform]
Let $\phi$ be a function defined on the space $\mathscr{P}^m$ of hyperplanes on $\R^m$. Then $\phi$ can be considered as an even function defined on the cylinder $\Sa^{m-1}\times \R$ which double-covers $\mathscr{P}^m$ by means of the map 
\[
(\m{w}, p) \in  \Sa^{m-1}\times \R \mapsto \{\m{x}\in\R^m: \langle\m{x},\m{w}\rangle=p\} \in \mathscr{P}^m.
\]
The dual Radon transform of $\phi$ on $\m{x}\in\R^{m}$ is defined as the average of $\phi$ on all hyperplanes passing through $\m{x}$, i.e.\
\[
\check{\phi}(\m{x}) = \frac{1}{\sigma_{m}} \int_{\Sa^{m-1}} \phi(\m{w}, \langle\m{x},\m{w}\rangle) \, dS_{\m{w}},
\]
where $\sigma_m=\frac{2\pi^{\frac{m}{2}}}{\Gam\left(\frac{m}{2}\right)}$ is the surface area of $\Sa^{m-1}$. This definition naturally extends to functions $f(x_0,\m{x})$ defined in $\R^{m+1}$. Indeed, making use of the double covering $\Sa^{m-1}\times \R$ of $\R^{m}\setminus \{0\}$ given by the transformation into spherical coordinates $\m{x}\mapsto r\m{w}$, we define the dual Radon transform of $f(x_0,\m{x})$ as
\[
\check{R}[f](x_0,\m{x}) = \frac{1}{\sigma_{m}} \int_{\Sa^{m-1}} f(x_0, \langle\m{x},\m{w}\rangle \m{w}) \, dS_{\m{w}}.
\]
\end{defi}

As mentioned before, in the work  \cite{MR3412341}, it was shown that the dual Radon transform $\check{R}$ maps slice monogenic functions into axial monogenic ones. Therefore, it makes sense to ask whether  $\check{R}$ is connected with $\tau_m$ and if so, how this connection takes place. The answer to these questions comes from our Theorem \ref{GenRel} and the following plane wave decomposition of the generalized CK-extension. Versions of this result has been proved  in \cite{MR3556034, MR3077647}. Nevertheless, we provide here a different proof which requires less computations and it is  easily adaptable to the most general version of the generalized CK-theorem, which deals with monogenic extensions of functions defined in subspaces of $\R^m$ of arbitrary co-dimensions, see Theorem 5.1.1 in \cite{MR1169463}.

\begin{thm}[Plane wave decomposition of the generalized CK-extension]\label{PWDGCK}
Let $\Om_1\inc \R$ be a real domain and let $f_0 \in \mathcal{A}(\Om_1)\otimes \R_m$. Then 
\begin{align*}
\textup{GCK}[f_0](x_0, \underline{x}) &= \frac{1}{\sigma_{m}} \left( \int_{\mathbb{S}^{m-1}} \exp(\langle \underline{\omega}, \underline{x} \rangle) \underline{\omega} \partial_{x_0}) dS_{\underline{\omega}} \right) f_0(x_0)=\frac{1}{\sigma_{m}} \int_{\mathbb{S}^{m-1}}f(x_{0}+ \langle \underline{x}, \underline{\omega} \rangle \underline{\omega}) dS_{\underline{\omega}},
\end{align*}
where $f = S[f_0]$ is the slice monogenic extension of $f_0$ to an axially symmetric slice neighbourhood  $\Om\inc\R^{m+1}$ of $\Om_1$, see (\ref{sliceExt}). 

Therefore the dual Radon transform provides a link between the operator $\textup{GCK}$ and the slice extension map $S$. This is $\textup{GCK} = \check{R} \circ S$, or equivalently,
\begin{equation}\label{ComRad}
\begin{tikzcd}[row sep = 2em, column sep = 8em]
 &\mathcal{AM}(\Om)  \arrow[dd, "\check{R}", leftarrow] \\
\mathcal{A}(\Om_1)\otimes \R_m \arrow[ur, "\textup{GCK}", rightarrow]   \arrow[dr, "S", rightarrow]  & \\
& \mathcal{SM}(\Om)
\end{tikzcd}.
\end{equation}
These three maps are right-module isomorphisms. In particular one has that $\check{R}= \textup{GCK} \circ S^{-1}$ and $\check{R}= S \circ \textup{GCK}^{-1}$.
\end{thm}
\begin{proof}
We begin by recalling the following particular cases of the Funk-Hecke Theorem (see e.g.\ \cite[Chp.~9.7]{MR1688958}) which are central for this proof. For suitable constants $C_0, C_1\in\R$ one has
\begin{align}\label{RadFrom}
\int_{\Sa^{m-1}} \langle {\m{x}}, {\m{\omega}} \rangle^j \, dS_{\m{\omega}} & = \begin{cases} C_0 |\m{x}|^j & j \mbox{ even} \\ 0, & j \mbox{ odd}\end{cases},  
& \int_{\Sa^{m-1}} \langle {\m{x}}, {\m{\omega}} \rangle^j \m{\omega} \, dS_{\m{\omega}} & = \begin{cases} 0 & j \mbox{ even} \\ C_1 |\m{x}|^{j-1}, & j \mbox{ odd}\end{cases}.
 \end{align}
The main idea of this proof is to construct a monogenic plane wave extension of the initial function $f_0 \in \mathcal{A}(\Om_1)\otimes \R_m$ and then transform it into an axial monogenic function via  a radialization of the form (\ref{RadFrom}). 
 
 \noindent To this end, let us consider the two-dimesional Cauchy-Riemann operator \[\pa_{\m{\omega}} = \pa_{x_0}+ \m{\omega}\langle \m{\omega}, \pa_{\m{x}}\rangle,\] where $\m{\omega}\in\Sa^{m-1}$ is a fixed direction in $\R^m$ and $\langle \m{\omega}, \pa_{\m{x}}\rangle= \sum_{j=1}^m \omega_j \pa_{x_j}$ is the directional derivative in the direction $\m{\omega}$. From (\ref{holExt}), it is known that the holomorphic extension of $f_0$ with respect to the operator $\pa_{\m{\omega}}$ is given by
\[
G(x_0, \langle \m{\omega}, {\m{x}}\rangle) = \exp\left(\langle \m{\omega}, {\m{x}}\rangle \m{\omega} \pa_{x_0}\right) f_0(x_0).
\]
This means that the identity $\pa_{\m{\omega}}G(x_0, \langle \m{\omega}, {\m{x}}\rangle) =0$ holds. 

\noindent On the other hand, the function $G(x_0, \langle \m{\omega}, {\m{x}}\rangle)$ is also in the kernel of the $(m+1)$-dimensional Cauchy-Riemann operator $\mathcal{D}_x=\pa_{x_0}+\pa_{\m{x}}$. Indeed, if we complete an orthonormal basis $\{\m{\omega}, \m{\xi}_1, \ldots, \m{\xi}_{m-1}\}$ of $\R^m$, then we can write $\mathcal{D}_x$ as 
\[
\mathcal{D}_x= \pa_{x_0}+\pa_{\m{x}} = \pa_{\m{\omega}}+ \sum_{j=1}^{m-1} \m{\xi}_j \langle \m{\xi}_j, \pa_{\m{x}}\rangle.
\]
The combination of this fact, with the identity
\[
\langle \m{\xi}_j, \pa_{\m{x}}\rangle G(x_0, \langle \m{\omega}, {\m{x}}\rangle) = \langle \m{\xi}_j, \m{\omega}\rangle \, \m{\omega} \pa_{x_0} \exp\left(\langle \m{\omega}, {\m{x}}\rangle \m{\omega} \pa_{x_0}\right) f_0(x_0) =0, 
\]
yields
\[
\mathcal{D}_x G(x_0, \langle \m{\omega}, {\m{x}}\rangle) = \pa_{\m{\omega}} G(x_0, \langle \m{\omega}, {\m{x}}\rangle) =0.
\]
From the Funk-Hecke formulas (\ref{RadFrom}), we easily obtain that the integral 
\[
\frac{1}{\sigma_{m}} \int_{\mathbb{S}^{m-1}}  G(x_0, \langle \m{\omega}, {\m{x}}\rangle)  dS_{\underline{\omega}} =\frac{1}{\sigma_{m}} \left( \int_{\mathbb{S}^{m-1}} \exp(\langle \underline{\omega}, \underline{x} \rangle) \underline{\omega} \partial_{x_0}) dS_{\underline{\omega}} \right) f_0(x_0),
\]
is a monogenic function of axial type whose restriction to the real line is given by $f_0(x_0)$. Finally, using the Taylor expansion of the holomorphic extension of $f_0$ around $x_0\in\R$ we obtain that
\begin{align*}
\textup{GCK}[f_0](x_0,\m{x})  &= \frac{1}{\sigma_{m}} \left( \int_{\mathbb{S}^{m-1}} \exp(\langle \underline{\omega}, \underline{x} \rangle) \underline{\omega} \partial_{x_0}) dS_{\underline{\omega}} \right) f_0(x_0)\\ &= 
\frac{1}{\sigma_{m}}  \int_{\mathbb{S}^{m-1}}\left( \sum_{j=0}^\infty \frac{\langle \m{\omega}, {\m{x}}\rangle^j \m{\omega}^j}{j!} \pa_{x_0}^j[f](x_0)\right) dS_{\underline{\omega}} \\
&= \frac{1}{\sigma_{m}} \int_{\mathbb{S}^{m-1}}f(x_{0}+ \langle \underline{x}, \underline{\omega} \rangle \underline{\omega}) dS_{\underline{\omega}},
\end{align*}
which completes the proof.
\end{proof}
Combining Theorem \ref{PWDGCK} with Theorem \ref{GenRel} we obtain the following plane wave decomposition of the Fueter-Sce-Qian mapping.
\begin{thm}\label{FullRelRadThm}
Let $ \Omega_2 \subset\C$ be an intrinsic complex domain and let $f: \Omega_2 \to \mathbb{C}$ be a holomorphic intrinsic function. Then for all dimensions $m \in \mathbb{N}$ we have
\begin{align*}
\tau_m \left[f(x_0+ \underline{x}) \right]&= \gam_m \; \textup{GCK}[f^{(m-1)}(x_{0})]\\
&= \frac{\gam_m }{\sigma_{m}} \int_{\mathbb{S}^{m-1}} f^{(m-1)}(x_{0}+ \underline{\omega} \langle \underline{x}, \underline{\omega} \rangle) dS_{\underline{\omega}},
\end{align*}
where $\gamma_m=   \frac{(-1)^{\frac{1-m}{2}} 2^{m-1} }{(m-1)!}  \Gamma \left( \frac{m+1}{2}\right)^2 $. In other words,  we have
\begin{equation}
\Delta^{\frac{m-1}{2}} \circ S= \gamma_m \, \textup{GCK} \circ \partial_{x_0}^{m-1}=\gamma_m \, \check{R} \circ S \circ \partial_{x_0}^{m-1},
\end{equation}
or equivalently,
\begin{equation}\label{DiagComp}
\begin{tikzcd}[row sep = 2em, column sep = 10em]
\mathcal{A}(\Om_1)\otimes \R_m \arrow[r, "S"] \arrow[rrd, "\tau_m \circ S", labels=below] \arrow[dd, "\gamma_m\; \pa_{x_0}^{m-1}", labels=left]&  \mathcal{SM}(\Om) \arrow[dr, "\tau_m"] &         \\
                                                                           &                                    &  \mathcal{AM}(\Om)                      \\
\mathcal{A}(\Om_1)\otimes \R_m \arrow[r, "S", labels=below] \arrow[rru, "\textup{GCK}"]  &  \mathcal{SM}(\Om)   \arrow[ur, "\check{R}", labels=below]             &      
\end{tikzcd}
\end{equation}
which completes the diagram (\ref{ComRad}).
\end{thm}

\begin{cor}[Plane wave decomposition of $Q^m_k$]
The elements of the family of monogenic monomials $\{Q^m_k\}_{k\in\N}$ defined in  (\ref{MalPol}) can be decomposed into plane waves as
\[
Q^m_k(x_0,\m{x}) = GCK[x_0^k] = \frac{1}{\sigma_{m}} \int_{\mathbb{S}^{m-1}} \left(x_{0}+ \langle \underline{x}, \underline{\omega} \rangle \underline{\omega}\right)^k dS_{\underline{\omega}}.
\]
\end{cor}

\begin{rem}
It is obvious that the similar plane waves decompositions can be obtained for the monogenic monomials of negative order $P^{(-k)}$ (see Definition \ref{DefMM}). Indeed, from formula (\ref{resu1}) we obtain
\begin{align*}
P^{(-k)} &= \frac{\lambda_m(m+k-2)!}{(k-1)! (m-1)!}\, \sgn(x_0)^{m+1} \, \textup{GCK}[x_0^{-k-m+1}] \\
&= \frac{\lambda_m(m+k-2)!}{(k-1)! (m-1)!}\, \frac{\sgn(x_0)^{m+1} }{\sigma_{m}}  \int_{\mathbb{S}^{m-1}} \left(x_{0}+ \langle \underline{x}, \underline{\omega} \rangle \underline{\omega}\right)^{-k-m+1} dS_{\underline{\omega}}.
\end{align*}
In particular, for $k=1$, the above formula yields the plane wave decomposition of the Cauchy kernel originally obtained in \cite{MR985370}. Indeed,
\begin{align*}
E(x) &= \frac{1}{\lan_m\sigma_{m+1}} P^{(-1)} = \frac{\sgn(x_0)^{m+1}}{\sigma_{m} \sigma_{m+1}}  \int_{\mathbb{S}^{m-1}} \left(x_{0}+ \langle \underline{x}, \underline{\omega} \rangle \underline{\omega}\right)^{-m} dS_{\underline{\omega}}.
\end{align*}
\end{rem}

\section{Coherent state transforms via the Fueter map}\label{CSTSect}
In \cite{MR3556034}, the authors introduced two extensions of the Segal-Bargmann or coherent state transform (CST) to the realm of Clifford analysis. These transforms result from extending the convolution with the heat kernel in $\R$ to $\R^{m+1}$ (instead of extending it to $\C$), by means of the slice extension map $S$ and the generalized CK-extension $\textup{GCK}$. In this way, one obtains two generalizations of the CST that map the right $\C_m$-module $\mathcal{L}^2(\R, dx_0)\otimes \C_m$ to the  modules of slice monogenic and axial monogenic functions in $\R^{m+1}$ respectively.  In a similar fashion to (\ref{ComRad}), these extensions were proven to be related by the Radon transform $\check{R}$. 

In this framework, one interesting question was left open in \cite[Remark 4.6]{MR3556034} regarding a possible alternative extension of the CST  by means of the Fueter-Sce-Qian map. In this section, we provide a complete answer to this question by defining such monogenic CST and studying its relation with those defined in \cite{MR3556034}. This will be possible thanks to diagram (\ref{DiagComp}), which offers a complete and comprehensive scheme on how the Fueter-Sce-Qian map, the slice monogenic extension, the generalized CK-extension, and the dual  Radon transform are connected.

Before we proceed, let us recall a few important notions about the Segal-Bargmann transform in $\R$, and its extensions to the Clifford analysis setting.

\subsection{Classical Segal-Bargmann transform in $\R$}{\phantom{a}}\newline
Let us consider first the heat equation in $\R$
\begin{equation}\label{HeatEq1}
\begin{cases}
\frac{1}{2}\Del_{0} [u] = \pa_{t} [u],\\[+.1cm]
u(0,x_0) = \phi(x_0),
\end{cases}
\end{equation}
where $\Del_{0}= \pa_{x_0}^2$ is the Euclidean Laplacian on the real line and $\phi$ is a ($\C$-valued) square integrable function on $\R$ with respect to Euclidean metric, i.e. $\phi\in\mathcal{L}^2(\R,dx_0)$. It is a well-known fact that the solution to this equation is given by 
\begin{equation}\label{HeatEq2}
u(t, x_0) = e^{\frac{t}{2}\Del_{0}} [\phi](x_0) = \int_{\R} \rho(t, x_0-y) \phi(y)\, d{y},
\end{equation}
where $\rho(t, x_0)=\frac{1}{(2\pi t)^{1/2}} e^{\frac{-x_0^2}{2t}}$ is the fundamental solution of (\ref{HeatEq1}), also known as the heat kernel in $\R$.

The Segal-Bargmann transform is defined as the convolution with the holomorphic extension of the heat kernel 
$\rho(1, x_0)$ to $\C$, i.e.\ it is the unique holomorphic extension of the convolution in (\ref{HeatEq2}).
\begin{defi}\label{SBDef}
The Segal-Bargmann or coherent state transform $U: \mathcal{L}^2(\R,dx_0) \fd Hol(\C)$ is defined as 
\[
U[f](z) = \int_{\R} \rho(1, z-y) f(y)\, d{y} = \frac{1}{(2\pi)^{1/2}}   \int_{\R}e^{\frac{-(z-y)^2}{2}}f(y)\, d{y}, \;\;\;\; z\in\C,
\]
where we recall that $Hol(\C)$ denotes the space of entire holomorphic functions in $\C$.
\end{defi}
It is well-known that the Segal-Bargmann transform $U$ is a unitary isomorphism from $\mathcal{L}^2(\R,dx_0)$ onto its image. Moreover, the image of $\mathcal{L}^2(\R,dx_0)$ under $U$ is known to be the subspace of $Hol(\C)$ composed of square integrable functions on $\C$ with respect to the measure $e^{-y^2}dx_0dy$. We denote this subspace as
\[
\mathcal{HL}^2 (\C, e^{-y^2}dx_0dy) = \left\{f(x_0+iy) \in\mathcal{H}(\C): \int_{\C} |f(x_0+iy)|^2 \; e^{-y^2}dx_0dy <\infty \right\}.
\]
From Definition \ref{SBDef}, it is easily seen that $U$ factorizes as follows 
\begin{equation}\label{ComSB}
\begin{tikzcd}[row sep = 2em, column sep = 8em]
                                                                                                                                                                       & \mathcal{HL}^2 (\C, e^{-y^2}dx_0dy) \arrow[d, "C", leftarrow]  \\
\mathcal{L}^2(\R,dx_0) \arrow[r, "e^{\frac{\Del_0}{2}}", rightarrow, labels=below] \arrow[ur, "U", rightarrow]  & \widetilde{\mathcal{A}}(\R)
\end{tikzcd}
\end{equation}
where $\widetilde{\mathcal{A}}(\R) \inc \mathcal{A}(\R)\otimes \C$ denotes the image of $\mathcal{L}^2(\R,dx_0)$ by the operator $e^{\frac{\Del_0}{2}}$, and $C$ denotes the holomorphic extension of elements $\widetilde{\mathcal{A}}(\R)$ to entire functions on $\C$, see (\ref{holExt}).

The slice monogenic and generalized CK extensions provide two natural ways of defining coherent state transforms in the Clifford setting by replacing the vertical arrow corresponding to the holomorphic extension $C$ in the diagram (\ref{ComSB}). This was the approach used in \cite{MR3556034}  to define the so-called slice monogenic and axial monogenic CSTs.  Let us briefly recall some important facts about these extensions.

\subsection{Slice monogenic and axial monogenic CSTs}{\phantom{a}}\newline
The main idea to define the slice monogenic CST (SCST), which we denote by $U_s$, is to replace the holomorphic extension $C$ from $\R$ to $\C$ in the vertical arrow of (\ref{ComSB}) by the slice extension map $S$, which leads to 
\begin{equation}\label{ComSM}
\begin{tikzcd}[row sep = 2em, column sep = 8em]
                                                                                                                                                                       & \mathcal{SM} (\R^{m+1})\otimes \C \arrow[d, "S", leftarrow]  \\
\mathcal{L}^2(\R,dx_0) \otimes \C_m \arrow[r, "e^{\frac{\Del_0}{2}}", rightarrow, labels=below] \arrow[ur, "U_s", rightarrow]  & \widetilde{\mathcal{A}}(\R)\otimes \C_m
\end{tikzcd}
\end{equation}
\begin{rem}\label{NeedForC1}
Although the action of $S$ was defined, in Section \ref{prem}, only on functions taking values on the real Clifford algebra $\R_m$, here we need to consider the action of $S$ on  $\C_m$-valued functions. Observe that this poses no problems in our argument since the definition of $S=\exp(\m{x}\pa_{x_0})$ can be naturally extended to $\C_m$-valued functions defined on $\R$. In this way, we obtain the map $S: \mathcal{A}(\R)\otimes \C_m \fd \mathcal{SM} (\R^{m+1})\otimes \C$.
\end{rem}
In \cite{MR3556034}, the explicit formula of the SCST was computed to be
\begin{align*}
U_s[f](x_0,\m{x}) &= \frac{1}{(2\pi)^{1/2}} \int_\R S\circ \exp\left(\frac{\pa_{x_0}^2}{2}\right) [e^{ipx_0}] \; \tilde{f}(p) dp \\
&=  \frac{1}{(2\pi)^{1/2}} \int_\R  e^{-\frac{p^2}{2}} \,  e^{ip(x_0+\m{x})} \; \tilde{f}(p) dp,
\end{align*}
where $\tilde{f}(p) = \frac{1}{(2\pi)^{1/2}} \int_\R e^{-ipx_0} f(x_0) \; dx_0$ is the Fourier transform of $f$ and 
\[
e^{ip(x_0+\m{x})} = e^{ipx_0} \left[\cosh(p|\m{x}|) + i\frac{x}{|\m{x}|} \sinh(p|\m{x}|)\right].
\]
We now summarize the main properties of this transform.
\begin{thm}[Thm 4.3 \cite{MR3556034}] \label{SMCSTthm}
Let $\mathcal{H}_s\inc \mathcal{SM}(\R^{m+1}) \otimes \C$ be the image of $\mathcal{L}^2(\R,dx_0)\otimes \C$ under $U_s$. Then the following statements hold.
\begin{itemize}
\item[$i)$] $\mathcal{H}_s \inc \mathcal{SML}^2(\R^{m+1}, dv_m)$ where  $\mathcal{SML}^2(\R^{m+1}, dv_m)$ is the Hilbert right module of functions in $\mathcal{SM}(\R^{m+1}) \otimes \C$ that are square integrable with respect to the measure $dv_m$ on $\R^{m+1}$ given by
\[
dv_m = \frac{2}{\sqrt{\pi}} \frac{1}{\sigma_m} \frac{e^{-|\m{x}|^2}}{|\m{x}|^{m-1}} dx_0 d{\m{x}},
\]
where $dx_0 d{\m{x}}= dx_0dx_1\cdots dx_m$ is the claasical Lebesgue measure in $\R^{m+1}$.
\item[$ii)$] The map $U_s$ in the diagram 
\begin{equation*}
\begin{tikzcd}[row sep = 2em, column sep = 8em]
                                                                                                                                                                       &\mathcal{H}_s \inc\mathcal{SM} (\R^{m+1})\otimes \C \arrow[d, "S", leftarrow]  \\
\mathcal{L}^2(\R,dx_0) \otimes \C_m \arrow[r, "e^{\frac{\Del_0}{2}}", rightarrow, labels=below] \arrow[ur, "U_s", rightarrow]  & \widetilde{\mathcal{A}}(\R)\otimes \C_m
\end{tikzcd}
\end{equation*}
is a unitary isomorphism for the inner product defined on $\mathcal{SML}^2(\R^{m+1}, dv_m)$ by means of the measure $dv_m$, i.e. for all $f,g\in \mathcal{L}^2(\R,dx_0) \otimes \C_m$ we have that
\[
\int_\R f(x_0)^\dagger g(x_0) \, dx_0 =  \frac{2}{\sqrt{\pi}} \frac{1}{\sigma_m} \int_{\R^{m+1}} U_s[f](x_0,\m{x})^\dagger \,  U_s[g](x_0,\m{x}) \, \frac{e^{-|\m{x}|^2}}{|\m{x}|^{m-1}} dx_0 d{\m{x}}. 
\]
where $\cdot^\dagger$ denotes the Hermitean conjugation (\ref{HermConj}) on $\C_m$.
\end{itemize}
\end{thm}

Similarly, the axial monogenic CST, which we denote by $U_a$, has been defined by replacing the vertical arrow in the diagram (\ref{ComSB}) by the generalized CK-extension. 
\begin{equation}\label{ComAM}
\begin{tikzcd}[row sep = 2em, column sep = 8em]
                                                                                                                                                                       & \mathcal{AM} (\R^{m+1})\otimes \C \arrow[d, "\textup{GCK}", leftarrow]  \\
\mathcal{L}^2(\R,dx_0) \otimes \C_m \arrow[r, "e^{\frac{\Del_0}{2}}", rightarrow, labels=below] \arrow[ur, "U_a", rightarrow]  & \widetilde{\mathcal{A}}(\R)\otimes \C_m
\end{tikzcd}
\end{equation}
\begin{rem}
In a similar fashion to Remark \ref{NeedForC1}, it is easily seen that the action of the operator GCK can be extended to $\C_m$-valued functions. In this way, we obtain GCK$:\mathcal{A}(\R)\otimes\C_m\fd \mathcal{AM} (\R^{m+1})\otimes \C$.
\end{rem}
Making use of the dual Radon transform and of the link it establishes between $\textup{GCK}$ and $S$, see (\ref{ComRad}), one obtains the following result.
\begin{thm}[Thm 4.5 \cite{MR3556034}]\label{Thm65}
The axial monogenic coherent state transform $U_a$ can be written as any of the following compositions
\[
U_a = \textup{GCK} \circ e^{\frac{\Del_0}{2}} = \check{R} \circ S \circ e^{\frac{\Del_0}{2}} = \check{R} \circ U_s.
\]
In particular, let $\mathcal{H}_a\inc \mathcal{AM}(\R^{m+1}) \otimes \C$ be the image of $\mathcal{L}^2(\R,dx_0)\otimes \C$ under $U_a$, then the restriction of the dual Radon transform $\check{R}$ defines an isomorphism between $\mathcal{H}_s$ and $\mathcal{H}_a$.  Moreover, the diagram
\[\begin{tikzcd}
	&&&&&& {\mathcal{H}_a \subset \mathcal{AM} (\R^{m+1})\otimes \mathbb{C}} \\
	{\mathcal{L}^2(\mathbb{R}, dx_0)\otimes\mathbb{C}_m} &&& {\widetilde{\mathcal{A}}(\mathbb{\R})\otimes\mathbb{C}_m} \\
	&&&&&& {\mathcal{H}_s \subset \mathcal{SM} (\R^{m+1})\otimes \mathbb{C}}
	\arrow["{e^{\frac{\Delta_0}{2}}}"{description}, from=2-1, to=2-4]
	\arrow["{\textup{GCK}}"', from=2-4, to=1-7]
	\arrow["{U_a}", from=2-1, to=1-7]
	\arrow["S", from=2-4, to=3-7]
	\arrow["{\check{R}}"', from=3-7, to=1-7]
	\arrow[from=2-1, to=3-7]
\end{tikzcd}\]

commutes and its exterior arrows are unitary isomorphisms, where the inner product on $\mathcal{H}_a$ is naturally induced by $\check{R}$ from the inner product on $\mathcal{H}_s$, i.e.\
\[
\langle F,G\rangle_{\mathcal{H}_a} = \int_{\R^{m+1}} \left(\check{R}\right)^{-1}\hspace{-.1cm}[F]^{\dagger}\, \left(\check{R}\right)^{-1}\hspace{-.1cm}[G] \;dv_m,
\]
where $dv_m$ is the measure on $\R^{m+1}$ defined in Theorem \ref{SMCSTthm}.
\end{thm}

\subsection{Axial CSTs via Fueter-Sce-Qian mapping theorem}{\phantom{a}}\newline
We are now in a position to introduce an axial CST by means of the Fueter-Sce-Qian mapping. Following the same principle as with the other extensions, this map would be defined by replacing the vertical arrow in (\ref{ComSB}) by the mapping $\tau_m\circ S$. which leads to the monogenic CST transform $\tau_m\circ S \circ e^{\frac{\Del_0}{2}}$. 

In order to uncover the relation of this transform with the CSTs $U_a$ and $U_s$, we first draw the following diagram which contains the mapping properties of all the CST's studied so far.
\[\begin{tikzcd}
	{\mathcal{L}^2(\mathbb{R}, dx_0)\otimes\mathbb{C}_m} &&& {\widetilde{\mathcal{A}}(\mathbb{\R})\otimes\mathbb{C}_m} &&& {\mathcal{H}_s \subset \mathcal{SM} (\R^{m+1})\otimes \mathbb{C}} \\
	\\
	&&&&&& {\mathcal{H}_a \subset \mathcal{AM} (\R^{m+1})\otimes \mathbb{C}} \\
	{\mathcal{L}^2(\mathbb{R}, dx_0)\otimes\mathbb{C}_m} &&& {\widetilde{\mathcal{A}}(\mathbb{\R})\otimes\mathbb{C}_m} \\
	&&&&&& {\mathcal{H}_s \subset \mathcal{SM} (\R^{m+1})\otimes \mathbb{C}}
	\arrow["{e^{\frac{\Delta_0}{2}}}"{description}, from=4-1, to=4-4]
	\arrow["{\textup{GCK}}"', from=4-4, to=3-7]
	\arrow["{U_a}", from=4-1, to=3-7]
	\arrow["S", from=4-4, to=5-7]
	\arrow["{\check{R}}"', from=5-7, to=3-7]
	\arrow[from=4-1, to=5-7]
	\arrow["{\tau_m}", from=1-7, to=3-7]
	\arrow["S", from=1-4, to=1-7]
	\arrow["{\tau_m\circ S}", from=1-4, to=3-7]
	\arrow["{e^{\frac{\Delta_0}{2}}}", from=1-1, to=1-4]
\end{tikzcd}\]

It is our goal to complete the above diagram and make it commutative, and to that end, our knowledge of the relation between $\tau_m\circ S$ and GCK will be crucial. Indeed, Theorem \ref{GenRel} clearly shows that the above diagram can be completed by adding the two missing vertical arrows (in the left and middle columns) corresponding to the operator $\gam_m \pa_{x_0}^{m-1}$. However, care must be taken since this leads to (weak) differentiation of $\mathcal{L}^2$-functions, which in general does not lead $\mathcal{L}^2$-functions as a result. This issue can be resolved if, instead of  $\mathcal{L}^2(\R,dx_0)$, we consider the Sobolev space $\mathcal{W}^{m-1,2}(\R, dx_0)$ in the right upper corner of diagram in Figure 1. We recall that Sobolev spaces are defined as
\[
\mathcal{W}^{k,2}(\R, dx_0) = \left\{f\in \mathcal{L}^2(\R,dx_0): \pa_{x_0}^{\el} [f] \in \mathcal{L}^2(\R,dx_0), \;\el=0,\ldots, k\right\}, \;\;\; k\in\N,
\]
where the above derivatives are considered in the weak sense. 

To complete the diagram in Figure 1, we need to check first that the weak derivative $\pa_{x_0}$ commutes with the smoothing operator $e^{\frac{\Del_0}{2}}$ when acting on $\mathcal{W}^{m-1,2}(\R, dx_0)$.
\begin{lem}\label{IntDiagmLem}
On any Sobolev space $\mathcal{W}^{k,2}(\R, dx_0)$,  $k\in\N$, the operator identity $e^{\frac{\Del_0}{2}} \circ \pa_{x_0} =  \pa_{x_0} \circ e^{\frac{\Del_0}{2}}$ holds. As a direct consequence, the following diagram commutes
\begin{equation}\label{IntDiagm}
\begin{tikzcd}[row sep = 3em, column sep = 6em]
\mathcal{W}^{m-1,2}(\R, dx_0) \arrow[r, "e^{\frac{\Del_0}{2}}", rightarrow] \arrow[d, " \pa_{x_0}^{m-1}", labels=left] & \widetilde{\mathcal{A}}(\R)^{*} \arrow[d, " \pa_{x_0}^{m-1}"] \\
 \mathcal{L}^2(\R,dx_0) \arrow[r, "e^{\frac{\Del_0}{2}}", rightarrow]                 &\widetilde{\mathcal{A}}(\R)
\end{tikzcd}
\end{equation}
where $\widetilde{\mathcal{A}}(\R)^{*}$  is the image of $\mathcal{W}^{m-1,2}(\R, dx_0)$ by $e^{\frac{\Del_0}{2}}$.
\end{lem}
\begin{proof}
This property follows from direct computations. Indeed, we recall that given $f\in \mathcal{W}^{k,2}(\R, dx_0)$ then
\[
e^{\frac{\Del_0}{2}}[f](x_0) = \frac{1}{(2\pi)^{1/2}}   \int_{\R}e^{\frac{-(x_0-y)^2}{2}}f(y)\, d{y}.
\]
Thus
\begin{align*}
\pa_{x_0} e^{\frac{\Del_0}{2}}[f](x_0) &=  \frac{1}{(2\pi)^{1/2}}   \int_{\R} \pa_{x_0} \left[e^{\frac{-(x_0-y)^2}{2}} \right] f(y)\, d{y}  \\
&=  \frac{-1}{(2\pi)^{1/2}}   \int_{\R} \pa_{y} \left[e^{\frac{-(x_0-y)^2}{2}} \right] f(y)\, d{y} \\
&=  \frac{1}{(2\pi)^{1/2}}   \int_{\R} e^{\frac{-(x_0-y)^2}{2}} \, f' (y)\, d{y} \\
&= e^{\frac{\Del_0}{2}}  \left[f'\right](x_0), 
\end{align*}
which proves the result.
\end{proof}
Combining the diagrams in Figure 1 and (\ref{IntDiagm}) with Theorem \ref{GenRel} we obtain the following commutative diagram which provides a full picture on the relation of the CSTs $U_a$ and $U_s$ and the Fueter-Sce-Qian mapping $\tau_m\circ S$.
\[\begin{tikzcd}
	{\mathcal{W}^{m-1,2}(\mathbb{R}, dx_0)\otimes\mathbb{C}_m} &&& {\widetilde{\mathcal{A}}(\mathbb{\R})^*\otimes\mathbb{C}_m} &&& {\mathcal{H}_s \subset \mathcal{SM} (\R^{m+1})\otimes \mathbb{C}} \\
	\\
	&&&&&& {\mathcal{H}_a \subset \mathcal{AM} (\R^{m+1})\otimes \mathbb{C}} \\
	{\mathcal{L}^2(\mathbb{R}, dx_0)\otimes\mathbb{C}_m} &&& {\widetilde{\mathcal{A}}(\mathbb{\R})\otimes\mathbb{C}_m} \\
	&&&&&& {\mathcal{H}_s \subset \mathcal{SM} (\R^{m+1})\otimes \mathbb{C}}
	\arrow["{e^{\frac{\Delta_0}{2}}}"{description}, from=4-1, to=4-4]
	\arrow["{\textup{GCK}}"', from=4-4, to=3-7]
	\arrow["{U_a}", from=4-1, to=3-7]
	\arrow["S", from=4-4, to=5-7]
	\arrow["{\check{R}}"', from=5-7, to=3-7]
	\arrow[from=4-1, to=5-7]
	\arrow["{\tau_m}", from=1-7, to=3-7]
	\arrow["S", from=1-4, to=1-7]
	\arrow["{\tau_m\circ S}", from=1-4, to=3-7]
	\arrow["{e^{\frac{\Delta_0}{2}}}", from=1-1, to=1-4]
	\arrow["{\gamma_m \partial_{x_0}^{m-1}}"', from=1-4, to=4-4]
	\arrow["{\gamma_m \partial_{x_0}^{m-1}}"', from=1-1, to=4-1]
\end{tikzcd}\]

Summarizing, we have obtained a new axial monogenic CST $\tau_m \circ S\circ e^{\frac{\Del_0}{2}}$ by means of the Fueter-Sce-Qian mapping. The relation of this transform with the previously defined $U_s$ and $U_a$ are given in the following theorem, which answers the question raised in \cite[Remark 4.6]{MR3556034}.
\begin{thm}\label{CSTDiagComThm}
The Fueter-Sce-Qian CST $\tau_m \circ S\circ e^{\frac{\Del_0}{2}}$ satisfies the following properties
\begin{align*}
\tau_m\circ S\circ e^{\frac{\Del_0}{2}} &= \gamma_m \, \textup{GCK} \circ \partial_{x_0}^{m-1} \circ e^{\frac{\Del_0}{2}}, & \mbox{on } &\;\;  \mathcal{L}^2(\R,dx_0) \otimes \C_m, \\
\tau_m\circ S\circ e^{\frac{\Del_0}{2}} &= \gamma_m \, U_a \circ \partial_{x_0}^{m-1} =  \gamma_m \, \check{R} \circ U_s \circ \partial_{x_0}^{m-1}, & \mbox{on } &  \;\; \mathcal{W}^{m-1,2}(\R, dx_0)\otimes \C_m.
\end{align*}
\end{thm}
\begin{proof}
The first formula directly follows from Theorem \ref{GenRel}, while the second formula is a direct consequence of  Lemma \ref{IntDiagmLem} and Theorem \ref{Thm65}.
\end{proof}

\bibliographystyle{abbrv}

\end{document}